\documentclass[11pt,twoside,letterpaper]{article} 
\usepackage{amssymb,amsmath}
\usepackage{times,fancyhdr}
\usepackage[dvips]{graphicx}


\makeatletter
\def\cleardoublepage{\clearpage\if@twoside \ifodd\c@page\else%
    \hbox{}%
    \thispagestyle{empty}%
    \newpage%
    \if@twocolumn\hbox{}\newpage\fi\fi\fi}
\makeatother

\def\figurename{Figure}
\makeatletter
\renewcommand{\fnum@figure}[1]{\figurename~\thefigure.}
\makeatother

\def\tablename{Table}
\makeatletter
\renewcommand{\fnum@table}[1]{\tablename~\thetable.}
\makeatother

\newtheorem{theo}{Theorem}[section]
\newtheorem{pr}[theo]{Proposition}

\newtheorem{remk}[theo]{Remark}

\newenvironment{rem}{\begin{remk}\normalfont}{\ \rule{0.5em}{0.5em}\end{remk}}

\newcommand{\refp}[1]{(\ref{#1})}
\numberwithin{equation}{section}

\renewcommand{\a}{\alpha}
\renewcommand{\b}{\beta}
\newcommand{\G}{\Gamma}
\renewcommand{\d}{\delta}
\newcommand{\D}{\Delta}
\newcommand{\z}{\zeta}
\renewcommand{\i}{\iota}
\renewcommand{\j}{\jmath}
\newcommand{\s}{\sigma}
\newcommand{\vs}{\varsigma}

\newcommand{\E}{\mathbb{E}}
\newcommand{\N}{\mathbb{N}}
\renewcommand{\P}{\mathbb{P}}
\newcommand{\R}{\mathbb{R}}
\newcommand{\Z}{\mathbb{Z}}

\newcommand{\bEx}{\mathbf{E}_{\mathbf{x}}}
\newcommand{\bEy}{\mathbf{E}_{\mathbf{y}}}
\newcommand{\bG}{\mathbf{G}}
\newcommand{\bH}{\mathbf{H}}
\newcommand{\bIE}{\mathbf{I}_{_{\boldsymbol{\mathcal{E}}}}}
\newcommand{\bIEprime}{\mathbf{I}_{_{\boldsymbol{\mathcal{E}'}}}}
\newcommand{\bIr}{\mathbf{I}_r}
\newcommand{\bK}{\mathbf{K}}
\newcommand{\bL}{\mathbf{L}}

\newcommand{\bp}{\mathbf{p}}
\newcommand{\bPo}{\mathbf{P}_{\!\mathbf{0}}}
\newcommand{\bPx}{\mathbf{P}_{\!\mathbf{x}}}
\newcommand{\bPy}{\mathbf{P}_{\!\mathbf{y}}}
\newcommand{\bPz}{\mathbf{P}_{\!\mathbf{z}}}

\newcommand{\bq}{\mathbf{q}}
\newcommand{\bQi}{\mathbf{Q}^{(\i)}}
\newcommand{\bQiun}{\mathbf{Q}^{(\i_1)}}
\newcommand{\bQik}{\mathbf{Q}^{(\i_k)}}
\newcommand{\bQikun}{\mathbf{Q}^{(\i_{k+1})}}
\newcommand{\br}{\mathbf{r}}
\newcommand{\bS}{\mathbf{S}}
\newcommand{\bT}{\mathbf{T}}
\newcommand{\bx}{\mathbf{x}}
\newcommand{\by}{\mathbf{y}}
\newcommand{\bz}{\mathbf{z}}

\newcommand{\cE}{\boldsymbol{\mathcal{E}}}
\newcommand{\cR}{\mathcal{R}}
\newcommand{\ind}{1\hspace{-.27em}\mbox{\rm l}}
\newcommand{\petitta}{\tau_{\raisebox{-.2ex}{$\scriptscriptstyle a$}}}
\newcommand{\petittab}{\tau_{\raisebox{-.2ex}{$\scriptscriptstyle a,b$}}}
\newcommand{\petittabc}{\tau_{\raisebox{-.2ex}{$\scriptscriptstyle a,b,c$}}}
\newcommand{\petittb}{\tau_{\raisebox{-.2ex}{$\scriptscriptstyle b$}}}
\newcommand{\petittc}{\tau_{\raisebox{-.2ex}{$\scriptscriptstyle c$}}}
\newcommand{\petittun}{\tau_{\raisebox{-.2ex}{$\scriptscriptstyle 1$}}}
\newcommand{\petittmoinsun}{\tau_{\raisebox{-.2ex}{$\scriptscriptstyle -1$}}}
\newcommand{\pxyi}{\mathbf{p}_{\mathbf{x},\mathbf{y}}^{(\i)}}
\newcommand{\pxyibis}{p_{x,y}^{(\i)}}
\newcommand{\pzyij}{\mathbf{p}_{\mathbf{z},\mathbf{y}}^{(\i-\j)}}
\newcommand{\qxyj}{\mathbf{q}_{\mathbf{x},\mathbf{y}}^{(\j)}}
\newcommand{\qxzj}{\mathbf{q}_{\mathbf{x},\mathbf{z}}^{(\j)}}
\newcommand{\rxkn}{\mathbf{r}_{\mathbf{x}}^{(k,n)}}
\newcommand{\rxon}{\mathbf{r}_{\mathbf{x}}^{(0,n)}}
\newcommand{\rykjn}{\mathbf{r}_{\mathbf{y}}^{(k-1,n-\j)}}
\newcommand{\so}{\s_{_0}}
\newcommand{\ta}{\tau_{\raisebox{-.3ex}{$\scriptstyle a$}}}
\newcommand{\tab}{\tau_{\raisebox{-.3ex}{$\scriptstyle a,b$}}}
\newcommand{\tabc}{\tau_{\raisebox{-.3ex}{$\scriptstyle a,b,c$}}}
\newcommand{\tac}{\tau_{\raisebox{-.3ex}{$\scriptstyle a,c$}}}
\newcommand{\tadeux}{\tau_{_{\scriptstyle a_2}}}
\newcommand{\tai}{\tau_{_{\scriptstyle a_i}}}
\newcommand{\taiimoinsun}{\tau_{_{\scriptstyle a_{i-1},a_i}}}
\newcommand{\taimoinsun}{\tau_{_{\scriptstyle a_{i-1}}}}
\newcommand{\taimoinsuniiplusun}{\tau_{_{\scriptstyle a_{i-1},a_i,a_{i+1}}}}
\newcommand{\taimoinsuniplusun}{\tau_{_{\scriptstyle a_{i-1},a_{i+1}}}}
\newcommand{\taiplusun}{\tau_{_{\scriptstyle a_{i+1}}}}
\newcommand{\taiiplusun}{\tau_{_{\scriptstyle a_i,a_{i+1}}}}
\newcommand{\tar}{\tau_{_{\scriptstyle a_r}}}
\newcommand{\tarmoinsunr}{\tau_{_{\scriptstyle a_{r-1},a_r}}}
\newcommand{\taj}{\tau_{_{\scriptstyle a_j}}}
\newcommand{\tajjprime}{\tau_{_{\scriptstyle a_j,a_{j'}}}}
\newcommand{\tajprime}{\tau_{_{\scriptstyle a_{j'}}}}
\newcommand{\tajsecond}{\tau_{_{\scriptstyle a_{j''}}}}
\newcommand{\tajjprimejsecond}{\tau_{_{\scriptstyle a_j,a_{j'},a_{j''}}}}
\newcommand{\taun}{\tau_{_{\scriptstyle a_1}}}
\newcommand{\taundeux}{\tau_{_{\scriptstyle a_1,a_2}}}
\newcommand{\tb}{\tau_{\raisebox{-.3ex}{$\scriptstyle b$}}}
\newcommand{\tbc}{\tau_{\raisebox{-.3ex}{$\scriptstyle b,c$}}}
\newcommand{\tbL}{\tilde{\mathbf{L}}}
\newcommand{\tc}{\tau_{\raisebox{-.3ex}{$\scriptstyle c$}}}
\newcommand{\tE}{\boldsymbol{\tau}_{_{\!\boldsymbol{\mathcal{E}}}}}
\newcommand{\tL}{\tilde{L}}
\newcommand{\tR}{\tau_{_{\!\mathcal{R}}}}

\newcommand{\ttai}{\tilde{\tau}_{_{\scriptstyle a_i}}}
\newcommand{\tzero}{\tau_{_0}}
\newcommand{\unE}{\mathbf{1\hspace{-.28em}l}_{_{\boldsymbol{\mathcal{E}}}}}
\newcommand{\unEprime}{\mathbf{1\hspace{-.28em}l}_{_{\boldsymbol{\mathcal{E}'}}}}
\newcommand{\unr}{\mathbf{1\hspace{-.28em}l}_r}
\newcommand{\vai}{\mathrm{var}_{a_i}}
\newcommand{\varpixa}{\varpi_{x,a}}
\newcommand{\varpixy}{\varpi_{x,y}}
\newcommand{\vso}{\vs_{_0}}

\newcommand{\Ea}{\E_a}
\newcommand{\Eadeux}{\E_{a_2}}
\newcommand{\Eai}{\E_{a_i}}
\newcommand{\Eamoinsun}{\E_{a-1}}
\newcommand{\Eaplusun}{\E_{a+1}}
\newcommand{\Eaun}{\E_{a_1}}
\newcommand{\Eb}{\E_b}
\newcommand{\Ebplusun}{\E_{b+1}}
\newcommand{\Ebmoinsun}{\E_{b-1}}
\newcommand{\Eo}{\E_{_0}}
\newcommand{\Ex}{\E_x}
\newcommand{\Pa}{\P_{\!a}}
\newcommand{\Pai}{\P_{a_i}}
\newcommand{\Pak}{\P_{a_k}}
\newcommand{\Pamoinsun}{\P_{\!a-1}}
\newcommand{\Paplusun}{\P_{\!a+1}}
\newcommand{\Pb}{\P_{b}}
\newcommand{\Pbmoinsun}{\P_{b-1}}
\newcommand{\Pbplusun}{\P_{b+1}}
\newcommand{\Po}{\P_{_{\!0}}}
\newcommand{\Px}{\P_{\!x}}
\newcommand{\Pxj}{\P_{\!x_j}}

\newcommand{\Tn}{\bT_{_{n,\boldsymbol{\mathcal{E}}}}}
\newcommand{\Tnprime}{\bT'_{\raisebox{0ex}{$\scriptscriptstyle n,\boldsymbol{\mathcal{E}}'$}}}
\newcommand{\Tnoprime}{\bT'_{\raisebox{0ex}{$\scriptscriptstyle n,0$}}}
\newcommand{\Tna}{T_{_{n,a}}}
\newcommand{\Tni}{\bT_{_{n-\i,\boldsymbol{\mathcal{E}}}}}
\newcommand{\Tnj}{\bT_{_{n-\j,\boldsymbol{\mathcal{E}}}}}
\newcommand{\Tno}{T_{_{n,0}}}
\newcommand{\TnR}{T_{_{n,\mathcal{R}}}}
\newcommand{\Tnuno}{T_{_{n-1,0}}}
\newcommand{\Tt}{\bT_{_{t,\boldsymbol{\mathcal{E}}}}}

\newcommand{\lqn}[1]{\noalign{\noindent $\displaystyle{#1}$}}

\begin{document}

\title{
{\begin{flushleft}
\vskip 0.45in
{\normalsize\bfseries\textit{Chapter~18}}
\end{flushleft}
\vskip 0.45in
\bfseries\scshape A random walk model related to the clustering of membrane receptors}}

\thispagestyle{fancy}
\fancyhead{}
\fancyhead[L]{In: Book Title \\
Editor: Editor Name, pp. {\thepage-\pageref{lastpage-01}}} 
\fancyhead[R]{ISBN 0000000000  \\
\copyright~2007 Nova Science Publishers, Inc.}
\fancyfoot{}
\renewcommand{\headrulewidth}{0pt}

\author{\bfseries\itshape Aim\'e LACHAL\thanks{E-mail address: aime.lachal@insa-lyon.fr}
\\
Universit\'e de Lyon, CNRS, INSA-Lyon, ICJ, UMR5208, F-69621, France
\\[1ex]
\footnotesize{Postal address:}\\[-.5ex]
\footnotesize{Institut National des Sciences Appliqu\'ees de Lyon}\\[-.5ex]
\footnotesize{P\^ole de Math\'ematiques/Institut Camille Jordan}\\[-.5ex]
\footnotesize{B\^atiment L\'eonard de Vinci, 20 avenue Albert Einstein}\\[-.5ex]
\footnotesize{69621 Villeurbanne Cedex, \textsc{France}}
}
\date{}
\maketitle
\thispagestyle{empty}
\setcounter{page}{1}

\label{lastpage-01}

\begin{abstract}
In a cellular medium, the plasmic membrane is a place
of interactions between the cell and its direct external environment.
A classic model describes it as a fluid mosaic.
The fluid phase of the membrane allows a lateral degree of freedom to
its constituents: they seem to be driven by random motions along the membrane.
On the other hand, experimentations bring to light inhomogeneities on the membrane;
these micro-domains (the so-called rafts) are very rich in proteins and phospholipids.
Nevertheless, few functional properties of these micro-domains have been
shown and it appears necessary to build appropriate models of the membrane
for recreating the biological mechanism.

In this article, we propose a random walk model
simulating the evolution of certain constituents--the so-called
ligands--along a heterogeneous membrane. Inhomogeneities--the rafts--are
described as being still clustered receptors. An important variable of interest to
biologists is the time that ligands and receptors bind during a fixed amount of time.
This stochastic time can be interpreted as a measurement of affinity/sentivity of ligands for receptors.
It corresponds to the sojourn time in a suitable set for a certain random walk.

We provide a method of calculation for the probability distribution
of this random variable and we next determine explicitly this distribution
in the simple case when we are dealing with only one ligand and one receptor.
We finally address some further more realistic models.
\end{abstract}



\vspace{.08in} \noindent {\textbf AMS Subject Classification:}
primary 60G50, 60J22;  secondary 60J10, 60E10.

\pagestyle{fancy}
\fancyhead{}
\fancyhead[EC]{Aim\'e Lachal}
\fancyhead[EL,OR]{\thepage}
\fancyhead[OC]{A random walk model related to the clustering of membrane receptors}
\fancyfoot{}
\renewcommand\headrulewidth{0.5pt}


\section{Introduction}

\subsection{The biological context}

In a cellular medium, the plasmic membrane is a place of interactions between
the cell and its direct external environment. A classic model describes it as
a fluid mosaic. The fluid phase of the membrane allows a lateral degree of
freedom to its constituents: they seem to be driven by random motions along
the membrane. As a first estimate, the membrane could
be viewed as a two-dimensional manifold on which the constituents are driven
by a Brownian motion. This makes the local concentration of constituents
independent from their localization on the membrane. In other words, the
membrane should be homogeneous. In fact, experimentations brought to light
inhomogeneities on the membrane. Indeed, different liquid phases were observed.
Certain constituents tend to group in clusters inside membrane domains,
forming dense receptor spots (sometimes called \textsl{rafts}) and depleted zones
elsewhere, instead of covering homogeneously the membrane surface.
These micro-domains are generally very rich in proteins and phospholipids,
but few functional properties have been shown.
Nowadays, the systematical presence of certain proteins and lipids
in different liquid phases has become a marker of such inhomogeneities.

Cellular response to changes in the concentration of different chemical species
(the so-called \textsl{ligands}) in the extracellular medium is induced
by the ligand binding to dedicated transmembrane receptors.
The receptor-ligand binding is based on local physical interactions,
ligand molecules randomly roam in the extracellular medium until they meet a receptor
at the cell surface and possibly dock. The binding mechanism, an important
feature for the biologists, provides a way of measurement of the
affinity/sensitivity of the ligands for the receptors.

We refer the reader to the paper by Car\'e \& Soula~\cite{hedi} for
an accurate description of the biological context.

In this article, we simplify the cellular medium by reducing it to
the dimension one and by discretizing Brownian motion, that is we propose a
naive one-dimensional random walk model to simulate the extracellular traffic
of ligands along a heterogeneous membrane. Inhomogeneities--the rafts--are
described as being a set of still monovalent receptors on the membrane.

An important variable of interest to biologists in the aforementioned binding
mechanism is the \textsl{docking-time}, this is the time that ligands and receptors
dock. Actually, this  variable is proportional to the time spent by ligands
on receptors during a fixed amount of time--the duration
of observation--with the rule that if several ligands
meet simultaneously the same receptor, they are counted only one time
because of the effect of the steric hindrance.

We propose a method of calculation for the probability distribution
of this random variable. Next, we shall consider particular cases and we shall
determine explicitly this distribution in the simple case when we are dealing
with only one ligand and one receptor.
Finally, we sketch more realistic models by introducing a circular random walk
or Brownian motion and we give some information about how to extend the results
we have obtained for the one-dimensional random walk to these cases.

\subsection{The random walk model}

The cellular boundary is modeled as the (discrete) integer line $\Z=\{\dots,-1,0,$ $1,\dots\}$.
Set $\N=\{0,1,2,\dots\}$ and $\N^*=\N\setminus\!\{0\}=\{1,2,\dots\}$.
We are given $r$~receptors $a_1,\dots,a_r \in\Z$ such that
$a_1<\dots<a_r$, $\cR=\{a_1,\dots,a_r\}$ and
$\ell$~ligands which are modeled as $\ell$~independent identically distributed
Bernoulli random walks with parameter $p\in(0,1)$.
The steps of each random walk take the value~$+1$ with probability~$p$ and
$-1$ with the probability~$q=1-p$.
The $\ell$~ligands thus induce an $\ell$-dimensional random walk $(\bS(\i))_{\i\in\N}$ on~$\Z^{\ell}$ with,
for any $\i\in\N$,
$$
\bS(\i)=(S_1(\i),\dots,S_\ell(\i)).
$$
Set $\displaystyle \cE=\bigcup_{j=1}^{\ell} \cE_{j}$ where the $\cE_j$'s are defined by
$$
\cE_{j}=\{(x_1,\dots,x_{\ell})\in\Z^{\ell} : x_j\in\cR\}=\Z^{j-1}\times\cR \times\Z^{\ell-j}.
$$
The $\cE_j$'s can be viewed as unions of hyperplanes parallel to the $j$-th
axis of coordinates; indeed,
$\displaystyle \cE_j=\bigcup_{i=1}^{r} \cE_{ij}$ where the $\cE_{ij}$'s
are the hyperplanes
$$
\cE_{ij}=\{(x_1,\dots,x_{\ell})\in\Z^{\ell} : x_j=a_i\}=\Z^{j-1}\times\{a_i\}\times\Z^{\ell-j}.
$$
The time spent by ligands on certain receptors during a fixed amount of time, say~$n$
($n\in\N^*$), is given by
$$
\Tn=\#\{\i\in\{1,\dots,n\}:\bS(\i)\in\cE\}=\sum_{\i=1}^n \ind_{\{\bS(\i)\in \cE\}},
$$
this is the sojourn time of the random walk $(\bS(\i))_{\i\in\N}$ in the set $\cE$
up to time~$n$. Indeed, $\bS(\i)\in \cE$ means that there exists two indices
$i\in\{1,\dots,r\}$ and $j\in\{1,\dots,\ell\}$ such that $S_j(\i)=a_i$, that is,
at least one ligand is located at one receptor. Moreover, if several ligands
$S_{j_1},\dots,S_{j_k}$ meet simultaneously a same receptor~$a_i$ at a certain
time~$\i$, i.e., $S_{j_1}(\i)=\dots=S_{j_k}(\i)=a_i$, then $\displaystyle \bS(\i)\in\bigcup_{j=1}^{\ell} \cE_{ij}$
and they are thus counted only one time in the quantity $\ind_{\{\bS(\i)\in \cE\}}$.

It is clear that the computations of the expectation and the variance of $\Tn$ are easy.
Computing the probability distribution of $\Tn$ is quite more complicated.
Our aim is to describe a possible way for computing this latter. We provide
matrix equations which may be solved by using numerical schemes.

\subsection{Settings}

We shall adopt the following convention for the settings: the roman letters
will be related to the dimension~$1$ (scalars) while the bold letters will be
related to the dimension~$\ell$ (vectors or matrices). We shall also put the space
variables ($x,y,z\in\Z,\bx,\by,\bz\in\Z^{\ell}$, etc.) in subscript and the
time variables ($\i,\j,k,n$, etc.) in superscript.
The letters $i,j,l,m$, etc. will be used as indices, e.g., for labeling the receptors.
Be aware of the difference between $\i$ and $i$, between $\j$ and $j$.
The letters $p,q,r,\bp,\bq,\br$, etc. will be used for defining certain probabilities,
the letters $G,H,K,\bG,\bH,\bK$, etc. for defining certain generating functions or matrices.
the letters $S,\bS$ for defining random walks,
the letters $n,\tau,T,\bT$ for defining certain times.

The settings $\bPx$ and $\bEx$ denote the probability and expectation
related to the \mbox{$\ell$-dimensional} random walk $(\bS(\i))_{\i\in\N}$ started
at a point $\bx\in\Z^{\ell}$ at time~$0$.
When the index of a one-dimensional random walk $(S_j(\i))_{\i\in\N}$ will not be used, we shall relabel
it as any generic one-dimensional random walk: $(S(\i))_{\i\in\N}$. The settings
$\Px$ and $\Ex$ denote the probability and expectation related to the
random walk $(S(\i))_{\i\in\N}$ started at a point $x\in\Z$ at time~$0$.

We introduce the first hitting time of the random walk $(\bS(\i))_{\i\in\N}$ in~$\cE$:
$$
\tE=\min\{\i\in \N^*: \bS(\i)\in \cE\}
$$
if there exists an index~$\i$ such that $\bS(\i)\in \cE$;
else we set $\tE=+\infty$. Similarly, we introduce the
first hitting times of the single or multiple levels $a_j$, $a_{j'}$, $a_{j''}$
together with that of the set~$\cR$ for the random walk $(S(\i))_{\i\in\N}$:
\begin{align*}
\taj
&=\min\{\i\in \N^*: S(\i)=a_j\},
\\
\tajjprime
&=\min\{\i\in \N^*: S(\i)\in\{a_j,a_{j'}\}\}=\min(\taj,\tajprime),
\\
\tajjprimejsecond
&=\min\{\i\in \N^*: S(\i)\in\{a_j,a_{j'},a_{j''}\}\}=\min(\taj,\tajprime,\tajsecond),
\\
\tR
&=\min\{\i\in \N^*: S(\i)\in\cR\}=\min(\taun,\dots,\tar)
\end{align*}
with the same convention: $\min(\emptyset)=+\infty$.

Let us define several family of probabilities:
for $\bx,\by\in\Z^{\ell}$, $\i,k,n\in\N$ and $\j\in\N^*$,
\begin{align*}
\pxyi
&=\bPx\{\bS(\i)=\by\}=\bPo\{\bS(\i)=\by-\bx\},
\\
\qxyj
&=\bPx\{\tE=\j,\bS(\tE)=\by\},
\\
\rxkn
&=\bPx\{\Tn=k\},
\end{align*}
We plainly have $\qxyj=0$ if $\by\notin\cE$ and
$\rxkn=0$ if $k>n$.
We also introduce their related generating functions: for $|u|,|v|<1$,
\begin{align*}
G_{\bx,\by}(u)
&=\sum_{\i=0}^{\infty} \pxyi u^{\i},
\\
H_{\bx,\by}(u)
&=\sum_{\j=1}^{\infty} \qxyj u^{\j}=\bEx(u^{\tE},\bS(\tE)=\by),
\\
K_{\bx}(u,v)
&=\sum_{k,n\in\N: \atop k\le n} \rxkn u^k v^n
=\sum_{n=0}^{\infty} \bEx(u^{\Tn} v^n),
\end{align*}
together with their associated matrices:
$$
\bG(u)=(G_{\bx,\by}(u))_{_{\scriptstyle\bx,\by\in \cE}},\quad
\bH(u)=(H_{\bx,\by}(u))_{_{\scriptstyle\bx,\by\in \cE}},\quad
\bK(u,v)=(K_{\bx}(u,v))_{_{\scriptstyle\bx\in \cE}}.
$$
The matrices $\bG$ and $\bH$ are infinite squared-matrices and $\bK$ is an
infinite column-matrix.

\subsection{Background on $1\mathrm{D}$-random walk}\label{background}

In this part, we supply several well-known results concerning classical
distributions related to the one-dimensional Bernoulli random walk.
We refer to~\cite{feller}, \cite{renyi} or~\cite{spitzer}.
In order to facilitate the reading and to make the paper self-contained,
we shall provide many details. In particular, we shall focus on the proofs which
use generating functions since they serve as a model of the main tool we
shall use throughout the paper. The proofs which do not involve any generating
functions will be postponed to Appendices~\ref{appendix1} and~\ref{appendix2}.

\subsubsection{Location of the random walk}\label{background1}

Set, for $x,y\in\N$ and $\i\in\N$,
$$
\pxyibis=\Px\{S(\i)=y\},\quad
G_{x,y}(u)=\sum_{\i=0}^{\infty} \pxyibis u^{\i}.
$$
We have
$$
\pxyibis=\binom{\i}{\frac{\i+x-y}{2}} p^{(\i+y-x)/2} q^{(\i+x-y)/2}
=\binom{\i}{\frac{\i+x-y}{2}} (pq)^{\i/2} \left(\frac pq\right)^{\!(y-x)/2}
$$
with the convention that $\binom{\i}{\a}=0$ for $\a\notin\N$ or $\a>\i$.
By putting $\varpixy=p$ if $x>y$, $\varpixy=q$ if $x<y$, $\varpixy=1$ if $x=y$,
we can rewrite $\pxyibis$ as
$$
\pxyibis=\binom{\i}{\frac{\i+|x-y|}{2}} (pq)^{\i/2} \left(\frac{pq}{\varpixy^2}\right)^{\!|x-y|/2}
$$
and next
$$
G_{x,y}(u)=\left(\frac{pq}{\varpixy^2}\right)^{\!|x-y|/2}
\sum_{\i\in\N:\,\i\ge |x-y|,\atop \i \text{ and } x-y \text{ with same parity}}
\binom{\i}{\frac{\i+|x-y|}{2}} (pq)^{\i/2} u^{\i}.
$$
By performing the change of index $\i\mapsto 2\i+|x-y|$ in the foregoing sum,
we derive
\begin{equation}\label{sumG}
G_{x,y}(u)=\left(\frac{pqu}{\varpixy}\right)^{\!|x-y|} \sum_{\i=0}^{\infty} \binom{2\i+|x-y|}{\i} (pqu^2)^{\i}.
\end{equation}
In order to simplify this last sum, we shall make use of the hypergeometric function
$F\left(\frac{l+1}{2},\frac{l+2}{2};l+1;\z\right)$
defined, with the usual notation $a_m=a(a+1)(a+2)\dots(a+m-1)$, by
$$
F\bigg(\frac{l+1}{2},\frac{l+2}{2};l+1;\z\bigg)
=\sum_{m=0}^{\infty} \frac{(\frac{l+1}{2})_{_{\scriptstyle m}}
(\frac{l+2}{2})_{_{\scriptstyle m}}}{(l+1)_m m!}\,\z^m
=\sum_{m=0}^{\infty} \binom{2m+l}{m} \!\left(\frac{\z}{4}\right)^{\!m}\!.
$$
By invoking Formula~15.1.14 of~\cite{abram}, p.~556, namely
$$
F\bigg(\frac{l+1}{2},\frac{l+2}{2};l+1;\z\bigg)
=\frac{2^l}{\left(1+\sqrt{1-\z}\,\right)^l\sqrt{1-\z}}
=\frac{2^l}{\sqrt{1-\z}}\!\left(\frac{1-\sqrt{1-\z}}{\z}\right)^{\!l}\!,
$$
we derive the relationship
\begin{equation}\label{sum}
\sum_{m=0}^{\infty} \binom{2m+l}{m}\z^m
=\frac{1}{\sqrt{1-4\z}}\left(\frac{1-\sqrt{1-4\z}}{2\z}\right)^{\!l}\!.
\end{equation}
Set
$$
A(u)=\sqrt{1-4pqu^2},\quad B^+(u)=\frac{1+A(u)}{2pu}=\frac{2qu}{1-A(u)},\quad B^-(u)=\frac{1-A(u)}{2pu}.
$$
We deduce from~\refp{sum} and~\refp{sumG} the following expression of
$G_{x,y}(u)$:
\begin{equation}\label{expressionG}
G_{x,y}(u)=\begin{cases}
\displaystyle\frac{[B^-(u)]^{x-y}}{A(u)} &\mbox{if } x>y,
\\[1ex]
\displaystyle\frac{1}{A(u)} &\mbox{if } x=y,
\\[2ex]
\displaystyle\frac{[B^+(u)]^{x-y}}{A(u)}
&\mbox{if } x<y.
\end{cases}
\end{equation}

\subsubsection{Hitting times of the random walk}\label{background2}

Let us now consider the family of hitting times related to $(S(\i))_{\i\in\N}$
started at $x\in\Z$: for any $a,b,c\in\Z$ such that $a<b<c$, set
\begin{align*}
\ta
&
=\min\{\i\ge 1: S(\i)=a\},
\\
\tab
&
=\min\{\i\ge 1: S(\i)\in\{a,b\}\},
\\
\tabc
&
=\min\{\i\ge 1: S(\i)\in\{a,b,c\}\}.
\end{align*}
We still adopt the convention $\min\emptyset =+\infty$.
In certain cases depending on the starting point~$x$, certain hitting times
are related. In fact, we have
$$
\tab=\begin{cases} \ta & \mbox{if $x<a$}, \\ \tb & \mbox{if $x>b$}, \end{cases}
$$
and
$$
\tabc=\begin{cases}
\ta & \mbox{if $x<a$,} \\
\tab & \mbox{if $x\in[a,b)$,} \\
\tbc & \mbox{if $x\in(b,c]$}, \\
\tc & \mbox{if $x>c$}.
\end{cases}
$$
Notice that when the starting point is $a$, time $\ta$ is the return time to level $a$;
when the starting point is $b$, time $\tabc$ depends on times $\tab$ and
$\tbc$. Let us introduce the probabilities
\begin{align*}
q_{x,a}^{(\j)}
&
=\Px\{\ta=\j\},
\\
q_{x,a,b}^{(\j)-}
&
=\Px\{\tab=\j,S(\tab)=a\}=\Px\{\ta=\j,\ta<\tb\},
\\
q_{x,a,b}^{(\j)+}
&
=\Px\{\tab=\j,S(\tab)=b\}=\Px\{\tb=\j,\tb<\ta\},
\end{align*}
together with their related generating functions:
\begin{align*}
H_{x,a}(u)
&
=\sum_{\j=1}^{\infty} q_{x,a}^{(\j)} u^{\j}=\Ex(u^{\petitta}),
\\
H_{x,a,b}^-(u)
&
=\sum_{\j=1}^{\infty} q_{x,a,b}^{(\j)-} u^{\j}=\Ex(u^{\petitta},\ta<\tb),
\\
H_{x,a,b}^+(u)
&
=\sum_{\j=1}^{\infty} q_{x,a,b}^{(\j)+} u^{\j}=\Ex(u^{\petittb},\tb<\ta).
\end{align*}

\noindent\textbf{\textsl{One-sided threshold}}
\vspace{.5\baselineskip}

For calculating the probability distribution of $\ta$, we invoke a ``continuity''
argument by observing that, for $\j\in\N^*$,
if $S(\j)=a$ then $\ta\le \j$, that is there exists an index $l\in\{1,\dots,\j\}$
such that $\ta=l$. This leads to the following relationship:
$$
\Px\{S(\j)=a\}=\Px\{S(\j)=a,\ta\le \j\}=\sum_{l=1}^{\j} \Px\{\ta=l\}\Pa\{S(\j-l)=a\}
$$
or, equivalently,
$$
p_{x,a}^{(\j)}=\sum_{l=1}^{\j} q_{x,a}^{(l)}p_{a,a}^{(\j-l)}.
$$
Then, the generating functions satisfies the equation
$$
G_{x,a}(u)=\d_{a,x}+H_{x,a}(u)G_{a,a}(u)
$$
from which we deduce, by~\refp{expressionG},
\begin{align}
H_{x,a}(u)
&=\begin{cases}
\displaystyle\frac{G_{x,a}(u)}{G_{a,a}(u)} &\mbox{if } x\ne a,
\\[3ex]
\displaystyle 1-\frac{1}{G_{a,a}(u)} &\mbox{if } x=a,
\end{cases}
\label{expressionH1}\\[2ex]
&
=\begin{cases}
[B^-(u)]^{x-a} &\mbox{if } x>a,
\\[2ex]
1-A(u) &\mbox{if } x=a,
\\[1ex]
[B^+(u)]^{x-a} &\mbox{if } x<a.
\end{cases}
\label{expressionH}
\end{align}
By using the hypergeometric function $F\left(\frac{l}{2},\frac{l+1}{2};l+1;\z\right)$
defined by
$$
F\bigg(\frac{l}{2},\frac{l+1}{2};l+1;\z\bigg)
=\sum_{m=0}^{\infty} \frac{(\frac{l}{2})_{_{\scriptstyle m}}
(\frac{l+1}{2})_{_{\scriptstyle m}}}{(l+1)_m m!}\,\z^m
=\sum_{m=0}^{\infty} \frac{l}{2m+l} \binom{2m+l}{m} \!\left(\frac{\z}{4}\right)^{\!m}\!.
$$
and referring to Formula~15.1.13 of~\cite{abram}, p.~556, namely
$$
F\bigg(\frac{l}{2},\frac{l+1}{2};l+1;\z\bigg)
=\frac{2^l}{\left(1+\sqrt{1-\z}\,\right)^l}
=2^l\left(\frac{1-\sqrt{1-\z}}{\z}\right)^{\!l}\!,
$$
we derive the relationship
\begin{equation}
\sum_{m=0}^{\infty} \frac{l}{2m+l} \binom{2m+l}{m}\z^m
=\left(\frac{1-\sqrt{1-4\z}}{2\z}\right)^{\!l}\!.
\label{sumbis}
\end{equation}
We then extract from~\refp{expressionH} and~\refp{sumbis}, for $x\ne a$,
\begin{align*}
H_{x,a}(u)
&
=\left(\frac{pqu}{\varpixa}\right)^{\!|x-a|} \sum_{\j=0}^{\infty}
\frac{|x-a|}{2\j+|x-a|}\binom{2\j+|x-a|}{\j} (pqu^2)^{\j}
\\
&
=\left(\frac{pq}{\varpixa^2}\right)^{\!|x-a|/2}
\sum_{\j\in\N:\,\j\ge |x-a|,\atop \j \text{ and } x-a \text{ with same parity}}
\frac{|x-a|}{\j}\binom{\j}{\frac{\j+|x-a|}{2}} (pq)^{\j/2} u^{\j}.
\end{align*}
We finally obtain
$$
q_{x,a}^{(\j)}=\frac{|x-a|}{\j}\,\binom{\j}{\frac{\j+|x-a|}{2}} (pq)^{\j/2}
\left(\frac{pq}{\varpixa^2}\right)^{\!|x-a|/2}
$$
that is, for $x\neq a$,
$$
q_{x,a}^{(\j)}=\frac{|x-a|}{\j}\,p_{x,a}^{(\j)}
=\frac{|x-a|}{\j}\,\binom{\j}{\frac{\j+x-a}{2}} p^{(\j+a-x)/2} q^{(\j+x-a)/2}.
$$
For $x=a$, we have, by considering the location of the first step of the walk,
$$
\Pa\{\ta=\j\}=p\,\Paplusun\{\ta=\j-1\}+q\,\Pamoinsun\{\ta=\j-1\}
$$
which leads to
$$
q_{a,a}^{(\j)}=\frac{1}{\j}\,\binom{\j}{\frac{\j+1}{2}} (pq)^{(\j+1)/2}.
$$

\noindent\textbf{\textsl{Two-sided threshold}}
\vspace{.5\baselineskip}

For calculating the probability distribution of $\tab$, we use
the strong Markov property by writing that for $x\in(a,b)$
\begin{align*}
H_{x,a}(u)
&
=\Ex(u^{\petitta})=\Ex(u^{\petitta},\ta<\tb)+\Ex(u^{\petitta},\tb<\ta)
\\
&
=\Ex(u^{\petitta},\ta<\tb)+\Ex(u^{\petittb},\tb<\ta)\,\Eb(u^{\petitta})
\\
&
=H_{x,a,b}^-(u)+H_{b,a}(u)\,H_{x,a,b}^+(u).
\end{align*}
Similarly,
$$
H_{x,b}(u)=H_{x,a,b}^+(u)+H_{a,b}(u)\,H_{x,a,b}^-(u)
$$
and we obtain that $H_{x,a,b}^+(u)$ and $H_{x,a,b}^-(u)$
solves the linear system
\begin{equation*}
\begin{cases}
H_{b,a}(u)\,H_{x,a,b}^+(u)+H_{x,a,b}^-(u)=H_{x,a}(u),
\\[1ex]
H_{x,a,b}^+(u)+H_{a,b}(u)\,H_{x,a,b}^-(u)=H_{x,b}(u),
\end{cases}
\end{equation*}
the solution of which writes
\begin{equation}\label{solutionH}
\begin{cases}
\displaystyle H_{x,a,b}^+(u)=\frac{H_{a,b}(u)\,H_{x,a}(u)-H_{x,b}(u)}{H_{a,b}(u)\,H_{b,a}(u)-1},
\\[3ex]
\displaystyle H_{x,a,b}^-(u)=\frac{H_{b,a}(u)\,H_{x,b}(u)-H_{x,a}(u)}{H_{a,b}(u)\,H_{b,a}(u)-1}.
\end{cases}
\end{equation}
Plugging the following expressions of $H_{x,a}(u)$ and $H_{x,b}(u)$
$$
\begin{cases}
H_{x,a}(u)=B^-(u)^{x-a} &\mbox{if } x>a
\\[1ex]
H_{x,b}(u)=B^+(u)^{x-b} &\mbox{if } x<b
\end{cases}
$$
into~\refp{solutionH}, we get that the generating function of $\tab$ is given by
$$
\Ex(u^{\petittab})=\Ex(u^{\petitta},\ta<\tb)+\Ex(u^{\petittb},\tb<\ta)
$$
where, for $x\in(a,b)$,
\begin{equation}\label{fg-tau-ab}
\begin{cases}
\displaystyle \Ex(u^{\petitta},\ta<\tb)=\frac{B^+(u)^{x-b}-B^-(u)^{x-b}}{B^+(u)^{a-b}-B^-(u)^{a-b}}
\\[2ex]
\displaystyle\Ex(u^{\petittb},\tb<\ta)=\frac{B^+(u)^{x-a}-B^-(u)^{x-a}}{B^+(u)^{b-a}-B^-(u)^{b-a}}
\end{cases}
\end{equation}
and then
\begin{equation}\label{fg-tau-ab-bis}
\Ex(u^{\petittab})=\frac{(1-B^-(u)^{b-a})B^+(u)^{x-a}-(1-B^+(u)^{b-a})B^-(u)^{x-a}}{B^+(u)^{b-a}-B^-(u)^{b-a}}.
\end{equation}
The generating functions~\refp{fg-tau-ab} can be inverted in order to write
out explicitly the distribution of $\tab$. The inversion hinges on the
decomposition of a certain rational fraction into partial fractions. The details
are postponed to Appendix~\ref{appendix1}. The result writes
\begin{align}
q_{x,a,b}^{(\j)-}
&
=2\left(\frac qp\right)^{\!(x-a)/2}
\Bigg[\sum_{l\in\N:\atop 1\le l<(b-a)/2} \cos^{x-a-1}\!\left(\frac{l\pi}{b-a}\right)
\sin\!\left(\frac{l\pi}{b-a}\right)
\nonumber\\
&
\hphantom{=\;}\sin\!\left(\frac{l(x-a)\pi}{b-a}\right)\!
\left(2\sqrt{pq}\cos\!\left(\frac{l\pi}{b-a}\right)\!\right)^{\j}\Bigg] u^{\j},
\nonumber\\
\label{distribution-tau-ab}
\\
q_{x,a,b}^{(\j)+}
&
=2\left(\frac pq\right)^{\!(x-b)/2}
\Bigg[\sum_{l\in\N:\atop 1\le l<(b-a)/2} \cos^{b-x-1}\!\left(\frac{l\pi}{b-a}\right)
\sin\!\left(\frac{l\pi}{b-a}\right)
\nonumber\\
&
\hphantom{=\;}\sin\!\left(\frac{l(b-x)\pi}{b-a}\right)\!
\left(2\sqrt{pq}\cos\!\left(\frac{l\pi}{b-a}\right)\!\right)^{\j}\Bigg] u^{\j}.
\nonumber
\end{align}

In the particular case $x=a$, we have
\begin{align*}
\begin{cases}
\Ea(u^{\petitta},\ta<\tb)
=pu\,\Eaplusun(u^{\petitta},\ta<\tb)+qu\,\Eamoinsun(u^{\petitta},\ta<\tb)
\\
\Ea(u^{\petittb},\tb<\ta)
=pu\,\Eaplusun(u^{\petittb},\tb<\ta)+qu\,\Eamoinsun(u^{\petittb},\tb<\ta)
\end{cases}
\end{align*}
which simplifies into
\begin{equation}\label{fg-tau-ab1}
\begin{cases}
\Ea(u^{\petitta},\ta<\tb)=pu\,\Eaplusun(u^{\petitta},\ta<\tb)+qu\,\Eo(u^{\petittun}),
\\
\Ea(u^{\petittb},\tb<\ta)=pu\,\Eaplusun(u^{\petittb},\tb<\ta).
\end{cases}
\end{equation}
Similarly, for $x=b$,
\begin{equation}\label{fg-tau-ab2}
\begin{cases}
\Eb(u^{\petitta},\ta<\tb)=qu\,\Ebmoinsun(u^{\petitta},\ta<\tb),
\\
\Eb(u^{\petittb},\tb<\ta)=qu\,\Ebmoinsun(u^{\petittb},\tb<\ta)+pu\,\Eo(u^{\petittmoinsun}).
\end{cases}
\end{equation}

\noindent\textbf{\textsl{Three-sided threshold}}
\vspace{.5\baselineskip}

At last, the probability distribution of $\tabc$ is characterized
by its generating function
\begin{align*}
\Ex(u^{\petittabc})
&
=\Ex(u^{\petitta},S(\tabc)=a)+\Ex(u^{\petittb},S(\tabc)=b)
\\
&
\hphantom{=\;}+\Ex(u^{\petittc},S(\tabc)=c)
\\
&
=\Ex(u^{\petitta},\ta<\tbc)+\Ex(u^{\petittb},\tb<\tac)+\Ex(u^{\petittc},\tau_c<\tab).
\end{align*}
We already observed that all the above generating functions
can be expressed by means of those of $\tab$ and $\tbc$ when
the starting point $x$ differs from $b$. When $x=b$, we have
\begin{align}
\Eb(u^{\petittabc},S(\tabc)=b)
&
=\Eb(u^{\petittb},\tb<\tac)
\nonumber\\
&
=pu\,\Ebplusun(u^{\petittb},\tb<\tau_c)+qu\,\Ebmoinsun(u^{\petittb},\tb<\ta).
\label{fg-tau-abc}
\end{align}
These probabilities can be expressed by means of~\refp{fg-tau-ab}.

\subsubsection{Stopped random walk}\label{background3}

In this part, we consider the families of generic ``stopping''--probabilities:
\begin{itemize}
\item
$\Px\{S(\j)=y,\j\le\ta\}$ ($a\in\R$) for $x,y\in(-\infty,a]$ or $x,y\in[a,+\infty)$;
\item
$\Px\{S(\j)=y,\j\le\tab\}$ ($a,b\in\R,a<b$) for $x,y\in[a,b]$;
\item
$\Px\{S(\j)=y,\j\le\tabc\}$ ($a,b,c\in\R,a<b<c$) for $x=y=b$.
\end{itemize}
These probabilities represent the distributions of the random walk stopped
when reaching level $a$, $a$ or $b$, $a$ or $b$ or $c$ respectively.
They can be evaluated by invoking the famous reflection principle.
Since this principle will not be used in our further analysis, we postpone
the intricate details to Appendix~\ref{appendix2}.
The result in the case of the one-sided barrier is, for $x,y<a$ or $x,y>a$,
\begin{align}
\Px\{S(\j)=y,\j\le \ta\}
&
=p_{x,y}^{(\j)}-\left(\frac pq \right)^{\!y-a} p_{x,2a-y}^{(\j)}
\nonumber\\
&
=\left[\!\binom{\j}{\frac{\j+x-y}{2}}
- \binom{\j}{\frac{\j+x+y}{2}-a}\!\right] p^{(\j+y-x)/2} q^{(\j+x-y)/2}.
\label{killed-a}
\end{align}
In the particular cases where $x=a$ or $y=a$, we have the following facts:
for $x\in\Z$ and $y=a$,
$$
\Px\{S(\j)=a,\j\le \ta\}=\Px\{\ta=\j\},
$$
and for $x=a$ and $y\in\Z$,
\begin{align*}
\Pa\{S(\j)=y,\j\le \ta\}
&
=p\,\Paplusun\{S(\j-1)=y,\j-1\le \ta\}
\\
&
\hphantom{=\;}
+q\,\Pamoinsun\{S(\j-1)=y,\j-1\le \ta\}.
\end{align*}
These last probabilities can be computed with the aid of~\refp{killed-a}.

In the case of the two-sided barrier, the result writes, for $x,y\in(a,b)$,
\begin{align}
\lqn{\Px\{S(\j)=y,\j\le \tab\}}
&
=\sum_{l=-\infty}^{\infty} \left(\frac pq \right)^{\!l(b-a)}
\left[p_{x,y-2l(b-a)}^{(\j)}-\left(\frac pq \right)^{\!y-a} p_{x,2a-2l(b-a)-y}^{(\j)}
\right]
\nonumber\\
&
=p^{(\j+y-x)/2} q^{(\j+x-y)/2}
\sum_{l=-\infty}^{\infty} \left[\!\binom{\j}{\frac{\j+x-y}{2}+l(b-a)}-
\binom{\j}{\frac{\j+x+y}{2}-a+l(b-a)}\!\right]\!.
\label{killed-ab}
\end{align}
The above sum actually is finite, limited to the indices $l$ such
that $\frac{y-x-\j}{2(b-a)} \le l\le \frac{y-x+\j}{2(b-a)}$ for the
first binomial coefficient and such that $\frac{2a-x-y-\j}{2(b-a)}
\le l\le \frac{2a-x-y+\j}{2(b-a)}$ for the second one.
In the particular case where $x\in\{a,b\}$ or $y\in\{a,b\}$, we have the following results.
For $x\in[a,b]$, if $y=a$,
$$
\Px\{S(\j)=a,\j\le \tab\}=\Px\{\ta=\j,\ta<\tb\}=q_{x,a,b}^{(\j)-},
$$
and, if $y=b$,
$$
\Px\{S(\j)=b,\j\le \tab\}=q_{x,a,b}^{(\j)+},
$$
These probabilities are given by~\refp{killed-a}. For $y\in[a,b]$, if $x=a$,
\begin{align*}
\Pa\{S(\j)=y,\j\le \tab\}
&
=p\,\Paplusun\{S(\j-1)=y,\j-1\le \tab\}
\\
&
\hphantom{=\;}
+q\,\Pamoinsun\{S(\j-1)=y,\j-1\le \ta\},
\end{align*}
and if $x=b$,
\begin{align*}
\Pb\{S(\j)=y,\j\le \tab\}
&
=p\,\Pbplusun\{S(\j-1)=y,\j-1\le \tb\}
\\
&
\hphantom{=\;}
+q\,\Pbmoinsun\{S(\j-1)=y,\j-1\le \tab\}.
\end{align*}
These probabilities are given by~\refp{killed-a} and~\refp{killed-ab}.

Finally, concerning the three-sided barrier, we simply have
\begin{align*}
\Pb\{S(\j)=b,\j\le \tabc\}
&
=p\,\Pbplusun\{S(\j-1)=b,\j-1\le \tbc\}
\\
&
\hphantom{=\;}
+q\,\Pbplusun\{S(\j-1)=b,\j-1\le \tab\}.
\end{align*}

\section{Methodology}

The aim of this part is to describe a method of calculation for the probability
distribution of $\Tn$ which could be numerically exploited.

\subsection{The probability distribution of $\bS(\i)$}

We have, for $\bx=(x_1,\dots,x_{\ell}),\by=(y_1,\dots,y_{\ell})\in \Z^{\ell}$,
\begin{align*}
\pxyi
&=\bPx\{\forall j\in\{1,\dots,\ell\},\,S_j(\i)=y_j\}
=\prod_{j=1}^{\ell} \Pxj\{S(\i)=y_j\}
\\
&=\Bigg[\prod_{j=1}^{\ell}\binom{\i}{\frac{\i+x_j-y_j}{2}}\!\Bigg]
p^{\left[\i+\sum_{j=1}^{\ell}(y_j-x_j)\right]/2}
q^{\left[\i+\sum_{j=1}^{\ell}(x_j-y_j)\right]/2}.
\end{align*}
In the above formula, $\pxyi$ does not vanish if and only if for
all $j\in\{1,\dots,\ell\}$, $x_j-y_j+\i$ is even and $|x_j-y_j|\le \i$.
Then the associated generating function is given by
\begin{equation}\label{G-inter}
G_{\bx,\by}(u)=\left(\frac pq\right)^{\left[\sum_{j=1}^{\ell}(y_j-x_j)\right]/2}
\sum_{\i=0}^{\infty}\Bigg[\prod_{j=1}^{\ell}\binom{\i}{\frac{\i+x_j-y_j}{2}}\!\Bigg]
\!\left(pqu^2\right)^{\i/2}.
\end{equation}
Set
$$
A(\xi_1,\dots,\xi_{\ell};z)
=\sum_{\i\in\N:\,\i\ge \max(|\xi_1|,\dots,|\xi_{\ell}|),
\atop \i,\xi_1,\dots,\xi_{\ell}\text{ with same parity}}
\Bigg[\prod_{j=1}^{\ell}\binom{\i}{\frac{\i+\xi_j}{2}}\!\Bigg]z^{\i/2}.
$$
The function $A$ does not vanish if and only if $\xi_1,\dots,\xi_{\ell}$ have
the same parity. By performing the change of index $\i\mapsto 2\i+|\xi_{j_{_0}}|$, where
$j_{_0}$ is an index such that $|\xi_{j_{_0}}|$ is the maximum of the
$|\xi_1|,\dots,|\xi_{\ell}|$, in the sum defining $A$, we get
$$
A(\xi_1,\dots,\xi_{\ell};z)
= z^{|\xi_{j_{_0}}|} \sum_{\i=0}^{\infty}
\Bigg[\prod_{j=1}^{\ell}\binom{2\i+|\xi_{j_{_0}}|}{\i+\frac{|\xi_{j_{_0}}|+\xi_j}{2}}\!\Bigg]z^{\i}.
$$
The quantity~\refp{G-inter} can be rewritten as follows.
%
\begin{pr}
For any $\bx=(x_1,\dots,x_{\ell}),\by=(y_1,\dots,y_{\ell})\in \Z^{\ell}$,
$$
G_{\bx,\by}(u)=\left(\frac pq\right)^{\left[\sum_{j=1}^{\ell}(y_j-x_j)\right]/2}
A\!\left(x_1-y_1,\dots,x_{\ell}-y_{\ell};pqu^2\right)\!.
$$
\end{pr}
%

We can express the function $A$ by means of hypergeometric functions.
To see this, we set $\a=|\xi_{j_{_0}}|$ and $\b_j=\frac{|\xi_{j_{_0}}|+\xi_j}{2}$
and assume, e.g., that $\xi_{j_{_0}}\le 0$ so that $\b_{j_{_0}}=0$.
In the case where $\xi_{j_{_0}}\ge 0$, we would have $\a-\b_{j_{_0}}=0$.
Invoking the duplication formula for the Gamma function, we write
\begin{align*}
\prod_{j=1}^{\ell}\binom{2\i+\a}{\i+\b_j}
&
=\prod_{j=1}^{\ell}\frac{\G(2\i+\a+1)}{\G(\i+\b_j+1)\G(\i+\a-\b_j+1)}
\\
&
=\frac{2^{(2\i+\a)\ell}}{\pi^{\ell/2}}
\prod_{j=1}^{\ell}\frac{\G(\i+\frac{\a+1}{2})\G(\i+\frac{\a+2}{2})}{\G(\i+\b_j+1)\G(\i+\a-\b_j+1)}
\\
&
=\frac{2^{(2\i+\a)\ell}}{\pi^{\ell/2}}
\frac{[\G(\i+\frac{\a+1}{2})\G(\i+\frac{\a+2}{2})]^{\ell}}{i!
\prod_{1\le j\le\ell,j\neq j_0}\G(\i+\b_j+1)
\prod_{1\le j\le\ell}\G(\i+\a-\b_j+1)}.
\end{align*}
Therefore, using the generalized hypergeometric function
\begin{align*}
\mbox{}_sF_t\Big(\!\begin{array}{c}\a_1,\dots,\a_s\\\b_1,\dots,\b_t\end{array}\!;z\Big)
&
=\sum_{m=0}^{\infty} \frac{(\a_1)_{_{\scriptstyle m}}\dots(\a_s)_{_{\scriptstyle m}}}{
(\b_1)_{_{\scriptstyle m}}\dots(\b_t)_{_{\scriptstyle m}}}\,\frac{z^m}{m!}
\\
&
=\frac{\G(\b_1)\dots\G(\b_t)}{\G(\a_1)\dots\G(\a_s)}
\sum_{m=0}^{\infty} \frac{\G(m+\a_1)\dots\G(m+\a_s)}{
\G(m+\b_1)\dots\G(m+\b_t)}\,\frac{z^m}{m!},
\end{align*}
we obtain
\begin{align*}
\lqn{A(\xi_1,\dots,\xi_{\ell};z)}
&
=\frac{2^{\a\ell}}{\pi^{\ell/2}}
\frac{[\G(\frac{\a+1}{2})\G(\frac{\a+2}{2})]^{\ell}}{
\prod_{1\le j\le\ell,j\neq j_0}\G(\b_j+1)
\prod_{1\le j\le\ell}\G(\a-\b_j+1)}\,z^{\a}
\\
&
\hphantom{=\;}\times\mbox{}_{2\ell}F_{2\ell-1}\!\left(\!\!\!\begin{array}{c}
(\a+1)/2,\dots,(\a+1)/2,(\a+2)/2,\dots,(\a+2)/2
\\
\b_1+1,\dots,\b_{\ell}+1,\a-\b_1+1,\dots,\a-\b_{\ell}+1\end{array}\!\!;4^{\ell} z\right)
\end{align*}
with the convention that in the list $(\a+1)/2,\dots,(\a+1)/2,(\a+2)/2,\dots,(\a+2)/2$
lying within the function $\mbox{}_{2\ell}F_{2\ell-1}$ above,the terms
$(\a+1)/2$ and $(\a+2)/2$ are repeated $\ell$~times and in the list
$\b_1+1,\dots,\b_{\ell}+1$ the term $\b_{j_{_0}}+1$ which equals one
is evicted.
Observing that the coefficient lying before the hypergeometric function
can be simplified into
\begin{align*}
\lqn{\frac{2^{\a\ell}}{\pi^{\ell/2}}
\frac{[\G(\frac{\a+1}{2})\G(\frac{\a+2}{2})]^{\ell}}{
\prod_{1\le j\le\ell,j\neq j_0}\G(\b_j+1)
\prod_{1\le j\le\ell}\G(\a-\b_j+1)}}
&
=\frac{\a^{\ell}}{\prod_{j=1}^{\ell}\b_j! \prod_{j=1}^{\ell}(\a-\b_j)!}
=\prod_{j=1}^{\ell} \binom{\a}{\b_j},
\end{align*}
we finally derive the following expression of $A(\xi_1,\dots,\xi_{\ell};z)$.
%
\begin{pr}
We have
\begin{align*}
\lqn{A(\xi_1,\dots,\xi_{\ell};z)=\Bigg(\prod_{j=1}^{\ell}
\binom{|\xi_{j_{_0}}|}{\b_j}\!\!\Bigg) z^{|\xi_{j_{_0}}|}}
\\
&
\hphantom{=\;}\times\!\mbox{}_{2\ell}F_{2\ell-1}\!\left(
\!\!\begin{array}{c} \frac{|\xi_{j_{_0}}|+1}{2}, \dots,\frac{|\xi_{j_{_0}}|+1}{2},
\frac{|\xi_{j_{_0}}|+2}{2},\dots,\frac{|\xi_{j_{_0}}|+2}{2}
\\[1ex]
\frac{|\xi_{j_{_0}}|+\xi_1}{2}+1,\dots,\frac{|\xi_{j_{_0}}|+\xi_{\ell}}{2}+1,
\frac{|\xi_{j_{_0}}|-\xi_1}{2}+1,\dots,\frac{|\xi_{j_{_0}}|-\xi_{\ell}}{2}+1
\end{array}\!\!;4^{\ell} z\right)\!.
\end{align*}
In the list $\frac{|\xi_{j_{_0}}|+\xi_1}{2}+1,\dots,\frac{|\xi_{j_{_0}}|+\xi_{\ell}}{2}+1,
\frac{|\xi_{j_{_0}}|-\xi_1}{2}+1,\dots,\frac{|\xi_{j_{_0}}|-\xi_{\ell}}{2}+1$
lying within the function $\mbox{}_{2\ell}F_{2\ell-1}$ above,
that of the two term $\frac{|\xi_{j_{_0}}|+\xi_{j_{_0}}}{2}+1$
and $\frac{|\xi_{j_{_0}}|-\xi_{j_{_0}}}{2}+1$ which equals one is evicted.
\end{pr}

\subsection{The probability distribution of $\tE$}

Fix $\j\in\N^*$. Notice that for $\by=(y_1,\dots,y_{\ell})\in \cE$, the event
$\{\tE=\j,\bS(\tE)=\by\}$ means that $\bS(\i)\notin \cE$ for all
$\i\in\{0,1,\dots,\j-1\}$, and $\bS(\j)=\by$.
Moreover, the event $\{\bS(\i)\notin\cE\}$ is equal to
$\displaystyle \bigcap_{j=1}^{\ell} \{S_j(\i)\notin\cR\}$ which means that
$S_j(\i)\notin\cR$ for all $j\in\{1,\dots,\ell\}$. Thus
\begin{align*}
\qxyj
&
=\bPx\{\forall \i\in\{0,1,\dots,\j-1\}, \forall j\in\{1,\dots,\ell\},\,S_j(\i)\notin \cR
\mbox{ and } S_j(\j)=y_j\}
\\
&
=\prod_{j=1}^{\ell} \Pxj\{\forall \i\in\{0,1,\dots,\j-1\},\,S_j(\i)\notin \cR
\mbox{ and } S_j(\j)=y_j\}
\\
&
=\prod_{j=1}^{\ell} \Pxj\{S(\j)=y_j,\j\le\tR\}.
\end{align*}
The quantity $\Pxj\{S(\j)=y_j,\j\le\tR\}$ is nothing but
the probability of the one-dimensional random walk stopped when reaching the set~$\cR$.
Notice that one of the $y_j$, $1\le j\le \ell$, at least, lies in $\cR$.
It can be depicted more precisely as follows.
Since the steps of the random walk are~$\pm1$, the random variable~$\tR$ under
the probability~$\Pxj$ is the first hitting time of the nearest neighbors
of~$x_j$ lying in~$\cR$, that is,
\begin{itemize}
\item
if $x_j\in(-\infty,a_1)$ (resp. $x_j\in(a_r,+\infty)$), then $\tR=\taun$ (resp. $\tR=\tar$);
\item
if $x_j\in(a_i,a_{i+1})$ for a certain index $i\in\{2,\dots,r-1\}$, then $\tR=\taiiplusun$;
\item
if $x_j=a_1$ (resp. $x_j=a_r$), then $\tR=\taundeux$ (resp. $\tarmoinsunr$);
\item
if $x_j=a_i$ for a certain index $i\in\{2,\dots,r-1\}$, then $\tR=\taimoinsuniiplusun$.
\end{itemize}
As a byproduct, we have
\begin{itemize}
\item
if $x_j\in(-\infty,a_1)$ and $y_j\in(-\infty,a_1]$ (resp. $x_j\in(a_r,+\infty)$
and $y_j\in[a_r,+\infty)$), then
\begin{align*}
&\Pxj\{S(\j)=y_j,\j\le\tR\}=\Pxj\{S(\j)=y_j,\j\le\taun\}
\\
\mbox{(resp. }&\Pxj\{S(\j)=y_j,\j\le\tR\}=\Pxj\{S(\j)=y_j,\j\le\tar\}\mbox{);}
\end{align*}
\item
if $x_j\in(a_i,a_{i+1})$ and $y_j\in[a_i,a_{i+1}]$ for a certain index $i\in\{1,\dots,r\}$,
$$
\Pxj\{S(\j)=y_j,\j\le\tR\}=\Pxj\{S(\j)=y_j,\j\le \taiiplusun\};
$$
\item
if $x_j=a_1$ and $y_j\in(-\infty,a_1)$ (resp. $x_j=a_r$ and $y_j\in(a_r,+\infty)$), then
\begin{align*}
&\Pxj\{S(\j)=y_j,\j\le\tR\}=\Pxj\{S(\j)=y_j,\j\le\taun\}
\\
\mbox{(resp. }&\Pxj\{S(\j)=y_j,\j\le\tR\}=\Pxj\{S(\j)=y_j,\j\le\tar\});
\end{align*}
\item
if $x_j=a_1$ and $y_j\in[a_1,a_2]$ (resp. $x_j=a_r$ and $y_j\in[a_{r-1},a_r]$), then
\begin{align*}
&\Pxj\{S(\j)=y_j,\j\le\tR\}=\Pxj\{S(\j)=y_j,\j\le\taundeux\}
\\
\mbox{(resp. }&\Pxj\{S(\j)=y_j,\j\le\tR\}=\Pxj\{S(\j)=y_j,\j\le\tarmoinsunr\});
\end{align*}
\item
if $x_j=a_i$ for a certain index $i\in\{2,\dots,r-1\}$ and $y\in[a_{i-1},a_i)$
(resp. $y\in(a_i,a_{i+1}]$), then
\begin{align*}
&\Pxj\{S(\j)=y_j,\j\le\tR\}=\Pxj\{S(\j)=y_j,\j\le\taiimoinsun\}
\\
\mbox{(resp. }&\Pxj\{S(\j)=y_j,\j\le\tR\}=\Pxj\{S(\j)=y_j,\j\le\taiiplusun\});
\end{align*}
\item
if $x_j=y_j=a_i$ for a certain index $i\in\{2,\dots,r-1\}$, then
$$
\Pxj\{S(\j)=y_j,\j\le\tR\}=\Pxj\{S(\j)=y_j,\j\le\taimoinsuniiplusun\}.
$$
\end{itemize}
In the other cases, the probability $\Pxj\{S(\j)=y_j,\j\le\tR\}$ vanishes.

All these probabilities belong to the following families of generic
``stopping''-probabilities which are explicitly given in Section~\ref{background3}.
With all this at hand, we can completely determine the joint
probability distribution of $(\tE,\bS(\tE))$.

The marginal probability distribution of $\tE$ can be easily related to that of $\tR$ according as
\begin{align*}
\bPx\{\tE\ge \j\}
&=\bPx\{\forall \i\in\{0,1,\dots,\j-1\}, \forall j\in\{1,\dots,\ell\},\,S_j(\i)\notin \cR\}
\\
&=\prod_{j=1}^{\ell} \Pxj\{\forall \i\in\{0,1,\dots,\j-1\},\,S_j(\i)\notin \cR\}
\\
&=\prod_{j=1}^{\ell} \Pxj\{\tR\ge \j\}.
\end{align*}
The probabilities $\Pxj\{\tR\ge \j\}$, $1\le j\le \ell$ can be easily
computed with the aid of the distributions described in Section~\ref{background2}.

We propose another possible way for deriving the distribution of $(\tE,\bS(\tE))$
which is characterized by the generating matrix~$\bH(u)$.
If $\bS(\i)\in\cE$, then $\tE\le \i$. Hence, using the strong Markov property,
we derive, for $\bx\in\Z^{\ell}$ and $\by\in\cE$, the relationship, for $\i\in\N^*$,
$$
\bPx\{\bS(\i)=\by\}=\sum_{\j=1}^{\i}\sum_{\bz\in \cE} \bPx\{\tE=\j,\bS(\tE)=\bz\} \bPz\{\bS(\i-\j)=\by\}
$$
or, equivalently,
$$
\pxyi=\sum_{\j=1}^{\i}\sum_{\bz\in \cE} \qxzj \pzyij.
$$
Therefore, taking the generating functions, for $\bx\in\Z^{\ell}$ and $\by\in\cE$,
\begin{align*}
G_{\bx,\by}(u)
&
=\d_{\bx,\by}+\sum_{\i=1}^{\infty} \pxyi u^{\i}
=\d_{\bx,\by}+ \sum_{\i,\j\in\N^*,\bz\in\cE:\atop \j\le \i} \qxzj \pzyij u^{\i}
\end{align*}
which can be rewritten as
\begin{equation}\label{eq-H}
G_{\bx,\by}(u)=\d_{\bx,\by}+\sum_{\bz\in\cE} H_{\bx,\bz}(u) G_{\bz,\by}(u),
\end{equation}
which leads, when restricting ourselves to $\bx\in\cE$,
to the matrix equation~\refp{eq-H-mat} below. Setting $\bIE$ for the identity matrix
$(\d_{\bx,\by})_{_{\scriptstyle \bx,\by\in\cE}}$, we have the following result.
%
\begin{theo}\label{prop-eq-H-mat}
The generating squared-matrix~$\bH(u)$ of the numbers $\qxyj $,
$\bx,\by\in\cE$, $\j\in\N$, which characterizes the joint probability distribution
of $(\tE,\bS(\tE))$, is a solution of the following matrix equation:
\begin{equation}\label{eq-H-mat}
[\bIE-\bH(u)]\bG(u)=\bIE.
\end{equation}
\end{theo}
%
This means that the generating functions $H_{\bx,\by}(u)$, $\bx,\by\in\cE$, are the solutions
of a system of an infinity of equations with an infinity of unknowns which seems
difficult to solve.

\subsection{The two first moments of $\Tn$}

We already observed that $\displaystyle\{\bS(\i)\notin\cE\}=\bigcap_{j=1}^{\ell}
\{S_j(\i)\notin\cR\}$. As a byproduct,
for $\bx=(x_1,\dots,x_{\ell})\in\Z^{\ell}$,
$$
\bPx\{\bS(\i)\notin \cE\}=\prod_{j=1}^{\ell}\Pxj\{S(\i)\notin \cR\}
$$
or
$$
\bPx\{\bS(\i)\in \cE\}=1- \prod_{j=1}^{\ell}(1-\Pxj\{S(\i)\in \cR\})
=1-\prod_{j=1}^{\ell} \bigg(1-\sum_{i=1}^r\Pxj\{S(\i)=a_i\}\bigg)\!.
$$
Now, the expectation of $\Tn$ can be easily computed as follows:
\begin{equation}\label{expectation}
\bEx(\Tn)=\sum_{\i=1}^n \bPx\{\bS(\i)\in \cE\}
=n-\sum_{\i=1}^n \prod_{j=1}^{\ell}\bigg(1-\sum_{i=1}^r\Pxj\{S(\i)=a_i\}\bigg)\!.
\end{equation}
The second moment of $\Tn$ could be evaluated as follows:
\begin{align*}
\bEx(\Tn^2)
&
=\sum_{\i=1}^n\sum_{\j=1}^n \bPx\{\bS(\i),\bS(\j)\in \cE\}
\\
&
=\sum_{\i=1}^n \bPx\{\bS(\i)\in \cE\}
+2\sum_{1\le\i<\j\le n} \bPx\{\bS(\i),\bS(\j)\in \cE\}.
\end{align*}
The foregoing double sum can be computed according as
\begin{align*}
\sum_{1\le\i<\j\le n} \bPx\{\bS(\i),\bS(\j)\in \cE\}
&
=\sum_{1\le\i<\j\le n} \sum_{\by\in\cE}\bPx\{\bS(\i)= \by\} \bPy\{\bS(\j-\i)\in \cE\}
\\
&
=\sum_{\i=1}^{n-1}\sum_{\by\in\cE} \bPx\{\bS(\i)=\by\} \sum_{\j=1}^{n-\i}\bPy\{\bS(\j)\in \cE\}
\\&
=\sum_{\i=1}^{n-1}\sum_{\by\in\cE} \bPx\{\bS(\i)=\by\} \bEy(\Tni).
\end{align*}
Consequently,
$$
\bEx(\Tn^2)=\bEx(\Tn)+2\sum_{\i=1}^{n-1}\sum_{\by\in\cE}
\bPx\{\bS(\i)=\by\} \bEy(\Tni).
$$
where $\bEx(\Tn)$ and the $\bEy(\Tni)$, $1\le \i\le n-1$, are given
by~\refp{expectation}.

\subsection{The probability distribution of $\Tn$}

We now propose a way for computing the distribution of~$\Tn$ under~$\bPx$ which is
determined by the family of numbers $\rxkn$, $0\le k\le n$.

For $1\le k\le n$, if $\Tn=k$, then $\tE\le n$, say $\tE=\j$ for a certain
$\j\in\{1,\dots,n\}$. Moreover, the sojourn time in~$\cE$ up to~$\tE$ is only~$1$,
and that after~$\cE$ up to~$n$, which is identical in distribution to
$\Tnj$, equals $k-1$. Hence, using the strong Markov property,
we derive the relationship, for $1\le k\le n$,
$$
\bPx\{\Tn=k\}=\sum_{\j=1}^n\sum_{\by\in \cE} \bPx\{\tE=\j,\bS(\tE)=\by\} \bPy\{\Tnj=k-1\}
$$
or, equivalently, for $1\le k\le n$,
$$
\rxkn=\sum_{\j=1}^n\sum_{\by\in \cE} \qxyj \rykjn.
$$
Therefore, taking the generating functions, for $\bx\in\Z^{\ell}$,
\begin{align}
K_{\bx}(u,v)
&=\sum_{n=0}^{\infty} \rxon v^n +\sum_{k,n\in\N^*: \atop k\le n} \rxkn u^k v^n
\nonumber\\
&=\sum_{n=0}^{\infty} \rxon v^n
+\sum_{k,n,\j\in\N^*,\,\by\in\cE: \atop k\le n\mbox{ \tiny and }\j\le n+1-k}
\qxyj \rykjn u^k v^n.
\label{fonc-gen-K}
\end{align}
On the other hand,
$$
\rxon=\bPx\{\Tn=0\}=\bPx\{\tE>n\}
$$
and the corresponding generating function is
\begin{align}
\sum_{n=0}^{\infty} \rxon v^n
&
=\sum_{n=0}^{\infty}\left(1-\sum_{k=0}^n \bPx\{\tE=k\} \right)\!v^n
\nonumber\\
&
=\frac{1}{1-v}-\sum_{k=0}^{\infty} \bPx\{\tE=k\} \left(\sum_{n=k}^{\infty} v^n\right)
\nonumber\\
&
=\frac{1}{1-v}-\sum_{k=0}^{\infty} \frac{v^k}{1-v}\,\bPx\{\tE=k\}
=\frac{1-\bEx(v^{\tE})}{1-v}.
\label{fonc-gen-tau}
\end{align}
In view of~\refp{fonc-gen-tau}, \refp{fonc-gen-K} can be rewritten as
\begin{equation}\label{eq-K}
K_{\bx}(u,v)=\frac{1-\bEx(v^{\tE})}{1-v}+u\sum_{\by\in\cE} H_{\bx,\by}(v) K_{\by}(u,v),
\end{equation}
This leads, when restricting ourselves to $\bx\in\cE$, to the following result.
%
\begin{theo}\label{prop-eq-K-mat}
The generating column-matrix $\bK(u,v)$ of the family of numbers $\bPx\{\Tn=k\}$,
$\bx\in\cE$, $k,n\in\N$, is a solution of the following matrix equation:
\begin{equation}\label{eq-K-mat}
[\bIE-u \bH(v)]\,\bK(u,v)=\frac{1}{1-v} \,[\unE-\tilde{\bH}(v)].
\end{equation}
The above matrix $\unE$ is the column-matrix consisting of~$1$,
that is $(1)_{_{\scriptstyle \bx\in\cE}}$, $\bH(v)$ is given by~\refp{eq-H-mat}
and $\tilde{\bH}(v)$ is defined by
$$
\tilde{\bH}(v)=\left(\bEx(v^{\tE})\right)_{_{\scriptstyle \bx\in\cE}}
=\Bigg(\sum_{\by\in\cE} H_{\bx,\by}(v)\Bigg)_{_{\!\!\scriptstyle \bx\in\cE}}=\bH(v)\unE.
$$
\end{theo}
%
In other words, the generating matrix $K_{\bx}(u,v)$, $\bx\in\cE$, are the solutions
of a system of an infinity of equations with an infinity of unknowns which seems
difficult to solve.
%
\begin{rem}
A slightly simpler equation may be obtained by setting
$$
\tilde{\bK}(u,v)=\bK(u,v)-\frac{1}{u(1-v)}\unE.
$$
In fact, by~\refp{eq-K-mat},
\begin{align*}
[\bIE-u \bH(v)]\tilde{\bK}(u,v)
&
=[\bIE-u\bH(v)] \!\left[\bK(u,v)-\frac{1}{u(1-v)}\unE\right]
\\
&
=\frac{1}{1-v}[\unE-\tilde{\bH}(v)]-\frac{1}{u(1-v)}[\unE-u\tilde{\bH}(v)].
\end{align*}
Thus, the modified generating function $\tilde{\bK}(u,v)$ satisfies
the matrix equation
\begin{equation}\label{eq-tildeK-mat}
[\bIE-u \bH(v)]\tilde{\bK}(u,v)=\frac{u-1}{u(1-v)}\unE.
\end{equation}
\end{rem}
%

Finally, for a starting point $\bx\in\Z^{\ell}\setminus \cE$, $K_{\bx}(u,v)$
can be expressed by means of the $G_{\bx,\by}(u)$ and $K_{\by}(u,v)$, $\by\in\cE$.
Indeed, by~\refp{eq-H}, we have for $\bx\in\Z^{\ell}\setminus \cE$ and
$\by\in\cE$,
$$
G_{\bx,\by}(u)=\sum_{\bz\in\cE} H_{\bx,\bz}(u) G_{\bz,\by}(u),
$$
that is, by introducing the row-matrices
$\bG_{\bx}(u)=(G_{\bx,\by}(u))_{\by\in\cE}$ and
$\bH_{\bx}(u)=(H_{\bx,\by}(u))_{\by\in\cE}$, the matrix $\bH_{\bx}(u)$
solves the equation
$$
\bG_{\bx}(u)=\bH_{\bx}(u)\bG(u).
$$
Next~\refp{eq-K} yields for $\bx\in\Z^{\ell}\setminus \cE$
$$
K_{\bx}(u,v)=\frac{1-\bEx(v^{\tE})}{1-v}+u \bH_{\bx}(u)\bK(u,v).
$$
\section{Particular cases}

In this part, we focus on the particular cases where $\ell=1$ or $r=1$. When $\ell=1$,
we are dealing with one ligand which can meet several receptors while in the case $r=1$,
we are concerned by one receptor which can be reached by several ligands.

\subsection{Case $\ell=1$}

In this case, our model is a one-dimensional random walk model and we have
$\cE=\cR$. We adapt the general settings by putting
$$
\TnR=\sum_{\i=1}^n \ind_{\{S(\i)\in \cR\}},
$$
and by writing, for $1\le i,j\le r$, $\i\in\N$ and $\j\in\N^*$, the following probabilities:
\begin{align*}
p_{i,j}^{(\i)}&=\Pai\{S(\i)=a_j\},\\
q_{i,j}^{(\j)}&=\Pai\{\tR=\j,S(\tR)=a_j\}=\Pai\{\tR=\taj=\j\},\\
r_i^{(k,n)}&=\Pai\{\TnR=k\}.
\end{align*}
We also introduce the generating functions
\begin{align*}
G_{i,j}(u)&=\sum_{\i=0}^{\infty} p_{i,j}^{(\i)} u^{\i},\\
H_{i,j}(u)&=\sum_{\j=1}^{\infty} q_{i,j}^{(\j)} u^{\j}=\E_{a_i}(u^{\tR},S(\tR)=a_j),\\
K_i(u,v)&=\sum_{k,n\in\N:\atop k\le n} r_i^{(k,n)}u^kv^n
=\sum_{n=0}^{\infty} \E_{a_i}\!\big(u^{\TnR}\big)v^n,
\end{align*}
together with the related (finite) matrices
$$
\bG(u)=(G_{i,j}(u))_{_{\scriptstyle 1\le i,j\le r}},\;
\bH(u)=(H_{i,j}(u))_{_{\scriptstyle 1\le i,j\le r}},\;
\bK(u,v)=(K_i(u,v))_{_{\scriptstyle 1\le i\le r}}.
$$
Referring to Section~\ref{background}, we explicitly have
$$
p_{i,j}^{(\i)}=\binom{\i}{\frac{\i+a_i-a_j}{2}} p^{(\i+\a_j-\a_i)/2} q^{(\i+\a_i-\a_j)/2}.
$$
$$
G_{i,j}(u)=\begin{cases}
\displaystyle\frac{[B^-(u)]^{a_i-a_j}}{A(u)} &\mbox{if } a_i>a_j,
\\[1ex]
\displaystyle\frac{1}{A(u)} &\mbox{if } a_i=a_j,
\\[2ex]
\displaystyle\frac{[B^+(u)]^{a_i-a_j}}{A(u)} &\mbox{if } a_i<a_j.
\end{cases}
$$
The matrix $\bH(u)$ is three-diagonal. More precisely, we have
$H_{i,j}(u)=0$ for $|i-j|\ge 2$ and
\begin{align*}
H_{i,i}(u)
&
=\Eai\big(u^{\tai},\tai<\taimoinsuniplusun\big),
\\
H_{i,i+1}(u)
&
=\Eai\big(u^{\taiplusun},\taiplusun<\tai\big),
\\
H_{i,i-1}(u)
&
=\Eai\big(u^{\taimoinsun},\taimoinsun<\tai\big).
\end{align*}
These quantities are explicitly given by~\refp{fg-tau-ab1}, \refp{fg-tau-ab2}
and~\refp{fg-tau-abc} in Section~\ref{background} which explicitly contains
the matrix $\bH(u)$.

The functions $H_{i,j}(u)$ can be also obtained by equations~\refp{eq-H}
which read here
$$
G_{i,j}(u)=\d_{i,j}+\sum_{k=1}^r H_{i,k}(u)G_{k,j}(u).
$$
This can be rewritten in terms of matrices as in~\refp{eq-H-mat}, by introducing the
unit $r\mathsf{x}r$-matrix $\bIr=(\d_{i,j})_{_{\scriptstyle 1\le i,j\le r}}$,
$$
\bG(u)=\bIr+\bH(u)\bG(u)
$$
and then
\begin{equation}\label{exp-H-with-G}
\bH(u)=\bIr-\bG(u)^{-1}.
\end{equation}
%
\begin{rem}
An alternative representation of $\bH(u)$ can be obtained as follows.
Let us introduce the probabilities
$$
\rho_{i,j}^{(\j)}=\Pai\{\taj=\j\},\quad \tilde{\rho}_{i,j}^{(\j)}
=\begin{cases} \rho_{i,j}^{(\j)} &\mbox{if } i\ne j\\ \d_{\j,0} &\mbox{if } i=j\end{cases}
$$
and their generating functions
$$
L_{i,j}(u)=\sum_{\j=1}^{\infty} \rho_{i,j}^{(\j)} u^{\j}=\Eai(u^{\taj}),
\quad\tL_{i,j}(u)=\sum_{\j=1}^{\infty} \tilde{\rho}_{i,j}^{(\j)} u^\j
=\begin{cases}L_{i,j}(u) &\mbox{if } i\ne j\\ 1&\mbox{if } i=j\end{cases}
$$
together with their related matrices
$$
\bL(u)=(L_{i,j}(u))_{_{\scriptstyle 1\le i,j\le r}},\quad
\tbL(u)=(\tL_{i,j}(u))_{_{\scriptstyle 1\le i,j\le r}}.
$$
Actually, the $\tilde{\rho}_{i,j}^{(\j)}$'s are associated with time
$\ttai=\min\{\i\in\N:S(\i)=a_i\}$. Notice that $\tai$ and $\ttai$ coincide
when the starting point differs from~$a_i$; in the case where the walk starts at~$a_i$,
$\ttai=0$ while $\tai$ is the first return time at~$a_i$.
We observe that
\begin{align*}
\rho_{i,j}^{(\j)}
&
=\Pai\{\tR\le \j,\taj=\j\}
\\
&
=\sum_{\i=1}^{\j-1} \sum_{k=1}^r \Pai\{\tR=\i,S(\tR)=a_k\} \Pak\{\taj=\j-\i\}
\\
&
\hphantom{=\;}+\sum_{k=1}^r \Pai\{\tR=\j,S(\tR)=a_j\}
=\sum_{\i=1}^{\j-1} \sum_{k=1}^r q_{i,k}^{(\i)} \,\tilde{\rho}_{k,j}^{(\j-\i)}.
\end{align*}
This implies the following relationship for the corresponding generating functions:
$$
L_{i,j}(u)=\sum_{k=1}^r H_{i,k}(u)\tL_{k,j}(u)
$$
or, equivalently,
$$
\bL(u)=\bH(u)\tbL(u),
$$
from which we deduce
$$
\bH(u)=\bL(u)\tbL(u)^{-1}.
$$
As a check, we notice that $\tbL(u)-\bL(u)$ is a diagonal matrix:
$$
\tbL(u)-\bL(u)=\mathrm{diag}(1-L_{i,i}(u))_{_{\scriptstyle 1\le i\le r}}=
\mathrm{diag}(1/G_{i,i}(u))_{_{\scriptstyle 1\le i\le r}}
$$
and then, by~\refp{expressionH1}, $\bG(u)[\tbL(u)-\bL(u)]=\tbL(u)$ which entails
$\bL(u)\tbL(u)^{-1}=\bIr-\bG(u)^{-1}$.
\end{rem}
%
Finally, the functions $K_i(v)$, $1\le i\le r$, are given by equation~\refp{eq-K}
which reads here
$$
K_i(u,v)=\frac{1-\Eai(v^{\tR})}{1-v}+u\sum_{j=1}^r H_{i,j}(v) K_j(u,v),
$$
and the equivalent matrix equation~\refp{eq-K-mat} writes
$$
[\bIr-u \bH(v)]\bK(u,v)=\frac{1}{1-v} \,[\unr-\tilde{\bH}(v)].
$$
The above column-matrix $\tilde{\bH}(v)$ is defined by
$$
\tilde{\bH}(v)=\left(\Eai(v^{\tR})\right)_{_{\scriptstyle 1\le i\le r}}
=\bH(v)\unr.
$$
The solution is the finite column-matrix
\begin{align}
\bK(u,v)
&
=\frac{1}{1-v} \,[\bIr-u \bH(v)]^{-1}[\unr-\tilde{\bH}(v)]
\nonumber\\
&
=\frac{1}{1-v} \,[\bIr-u \bH(v)]^{-1}[\bIr-\bH(v)]\unr.
\label{expressionK-mat}
\end{align}
%
\begin{rem}
Referring to Equation~\refp{eq-tildeK-mat}, we see that the modified
generating function $\tilde{\bK}(u,v)=\bK(u,v)-\frac{1}{u(1-v)}\unr$
is given by
\begin{equation}\label{expression-tildeK-mat}
\tilde{\bK}(u,v)=\frac{u-1}{u(1-v)}\,[\bIr-u \bH(v)]^{-1}\unr.
\end{equation}
\end{rem}
%
\begin{rem}
By inserting the expression~\refp{exp-H-with-G} of $\bH(v)$ by means
of $\bG(v)$ into~\refp{expressionK-mat}, we get that
\begin{align*}
[\bIr-u \bH(v)]^{-1}[\bIr-\bH(v)]
&
=[(1-u)\bIr+u \bG(v)^{-1}]^{-1}[\bG(v)^{-1}]
\\
&
=[(1-u)\bG(v)+u \bIr]^{-1}.
\end{align*}
So, $\bK(u,v)$ can be expressed in terms of $\bG(v)$ as
$$
\bK(u,v)=\frac{1}{1-v}\,[(1-u)\bG(v)+u \bIr]^{-1}\unr.
$$
This representation is simpler than~\refp{expressionK-mat}.
Nonetheless, it is not tractable for inverting the generating
function $\bK(u,v)$ with respect to $u$.
\end{rem}
%
Expanding $[\bIr-u \bH(v)]^{-1}$ into $\sum_{k=0}^{\infty} u^k \bH(v)^k$,
we get by~\refp{expressionK-mat} (or~\refp{expression-tildeK-mat})
$$
\bK(u,v)=\frac{1}{1-v} \sum_{k=0}^{\infty} u^k [\bH(v)^k (\bIr-\bH(v)) \unr]
$$
from which we extract the following proposition.
%
\begin{pr}
The probability distribution of $\TnR$ satisfies, for any $k\in\N$,
\begin{equation}\label{law-TnR}
\sum_{n=k}^{\infty} \left((\Pai\{\TnR=k\})_{_{\scriptstyle 1\le i\le r}}\right) v^n
=\frac{1}{1-v} \,\bH(v)^k [\bIr-\bH(v)] \unr
\end{equation}
where each sides of the equality are column-matrices.
\end{pr}
%
We could go further in the computations: expanding $1/(1-v)$ into
$\sum_{n=0}^{\infty} v^n$, we obtain by~\refp{law-TnR} that
$$
\sum_{n=k}^{\infty} \left(\!\left(r_i^{(k,n)}\right)_{1\le i\le r}\right)v^n=
\left(\sum_{n=k}^{\infty} v^n\right) \bH(v)^k [\bIr-\bH(v)]\unr.
$$
Introducing the matrix $\bQi=\big(q_{i,j}^{(\i)}\big)_{1\le i,j\le r}$, we rewrite $\bH(v)$ as
$\sum_{\i=1}^{\infty} \bQi v^{\i}$. Then,
$$
\bH(v)^k=\sum_{\i=k}^{\infty}\left(\sum_{\i_1,\dots,\i_k\in\N^*:\atop\i_1+\dots+\i_k=\i}
\bQiun\dots\bQik\right)\! v^{\i}
$$
and
\newcommand{\ghost}{\vphantom{\sum_{\i_1\atop\i_1}}}
\begin{align*}
\lqn{\bH(v)^k[\bIr-\bH(v)]}
=\sum_{\i=k}^{\infty}\left(\ghost\right.
\sum_{\i_1,\dots,\i_k\in\N^*:\atop\i_1+\dots+\i_k=\i}\bQiun\dots\bQik
-\sum_{\i_1,\dots,\i_{k+1}\in\N^*:\atop\i_1+\dots+\i_{k+1}=\i}\bQiun\dots\bQikun
\left.\ghost\right) \!v^{\i}.
\end{align*}
Next,
\begin{align*}
\lqn{\left(\sum_{n=k}^{\infty} v^n\right) \bH(v)^k [\bIr-\bH(v)]\unr}
=\sum_{n=k}^{\infty} \left[\sum_{\i=k}^n\left(\ghost\right.\right.
\sum_{\i_1,\dots,\i_k\in\N^*:\atop\i_1+\dots+\i_k=\i}\bQiun\dots\bQik\unr
-\sum_{\i_1,\dots,\i_{k+1}\in\N^*:\atop\i_1+\dots+\i_{k+1}=\i}\bQiun\dots\bQikun\unr
\left.\left.\ghost\right)\!\right] \!v^{\i}.
\end{align*}
We finally deduce the result below.
%
\begin{theo}
The probability $\Pai\{\TnR=k\}$ is the $i$-th term of the column-matrix
$$
\sum_{\i_1,\dots,\i_k\in\N^*:\atop k\le\i_1+\dots+\i_k\le n}\bQiun\dots\bQik\unr
-\sum_{\i_1,\dots,\i_{k+1}\in\N^*:\atop k\le\i_1+\dots+\i_{k+1}\le n}\bQiun\dots\bQikun\unr.
$$
\end{theo}
%
The expectation of $\TnR$ can be derived by evaluating the derivative of its
generating function at $u=1$. Indeed,
\begin{align*}
\frac{\partial \bK}{\partial u}(1,v)
&=\left.\frac{1}{1-v} \,\bH(v) [\bIr-u \bH(v)]^{-2}[\bIr-\bH(v)]\unr\right|_{u=1}
\\
&=\frac{1}{1-v} \,\bH(v) [\bIr-\bH(v)]^{-1}\unr.
\end{align*}
Since
$$
\bH(v) [\bIr-\bH(v)]^{-1}=[\bIr-(\bIr-\bH(v))] [\bIr-\bH(v)]^{-1}
=[\bIr-\bH(v)]^{-1}-\bIr,
$$
we get the expectation of $\TnR$ when the random walk starts at a receptor at time~$0$:
$$
(\Eai(\TnR))_{_{\scriptstyle 1\le i\le r}}=\frac{1}{1-v} \left([\bIr-\bH(v)]^{-1}-\bIr\right)\unr.
$$
Next, the second moment of $\TnR$ can be extracted by computing the second derivative of its
generating function at $u=1$:
\begin{align*}
\frac{\partial^2 \bK}{\partial u^2}(1,v)
&=\left.\frac{2}{1-v} \,\bH(v)^2 [\bIr-u \bH(v)]^{-3}[\bIr-\bH(v)]\unr\right|_{u=1}
\\
&=\frac{2}{1-v} \,\bH(v)^2 [\bIr-\bH(v)]^{-2}\unr.
\end{align*}
Since
\begin{align*}
\lqn{\bH(v)^2 [\bIr-\bH(v)]^{-2}+\bH(v) [\bIr-\bH(v)]^{-1}}
&=\bH(v) [\bIr-\bH(v)]^{-2} [\bH(v)+ (\bIr-\bH(v))]
\\
&=[\bIr-(\bIr-\bH(v))] [\bIr-\bH(v)]^{-2}
\\
&=[\bIr-\bH(v)]^{-2}-[\bIr-\bH(v)]^{-1},
\end{align*}
we get the variance of $\TnR$:
$$
(\vai(\TnR))_{_{\scriptstyle 1\le i\le r}}=\frac{1}{1-v}
\left([\bIr-\bH(v)]^{-2}-[\bIr-\bH(v)]^{-1}\right)\unr.
$$

\subsection{Case $r=1$}

In this case, we have one receptor: $\cR=\{a\}$, and then
$\displaystyle \cE=\bigcup_{j=1}^{\ell} \cE_{j}$ where
$\cE_{j}$ is the hyperplane $\Z^{j-1}\times\{a\} \times\Z^{\ell-j}$.
Formulas~\refp{eq-H-mat} and~\refp{eq-K-mat} in Theorems~\ref{prop-eq-H-mat}
and~\ref{prop-eq-K-mat} do not simplify.


\subsection{Case $\ell=1$ and $r=2$}

We have in this case $\cE=\cR=\{a_1,a_2\}$. We work with the modified generating
function $\tilde{\bK}(u,v)$ given by~\refp{expression-tildeK-mat}. We have to invert the
$2\mathsf{x}2$-matrix $\mathbf{I}_2-u\bH(v)$. Here, $\bH(v)$ reads
$$
\bH(v)=\begin{pmatrix}
H_{1,1}(v) & H_{1,2}(v) \\ H_{2,1}(v) & H_{2,2}(v)
\end{pmatrix}
=\begin{pmatrix}
\Eaun(v^{\tau_{a_1}},\taun<\tadeux) & \Eaun(v^{\tau_{a_2}},\tadeux<\taun)
\\
\Eadeux(v^{\tau_{a_1}},\taun<\tadeux) & \Eadeux(v^{\tau_{a_2}},\tadeux<\taun)
\end{pmatrix}.
$$
The entries of this matrix are given by~\refp{fg-tau-ab1} and~\refp{fg-tau-ab2}.
Setting
$$
\D(u,v)=1-u [H_{1,1}(v)+H_{2,2}(v)]+u^2[H_{1,1}(v)H_{2,2}(v)-H_{1,2}(v)H_{2,1}(v)],
$$
we have
$$
[\mathbf{I}_2-u\bH(v)]^{-1}=\frac{1}{\D(u,v)}
\begin{pmatrix}
1-u H_{2,2}(v) & u H_{2,1}(v)\\ u H_{1,2}(v)   & 1-u H_{1,1}(v)
\end{pmatrix}
$$
and then, by~\refp{expression-tildeK-mat},
$$
\tilde{\bK}(u,v)=\frac{u-1}{u(1-v)\D(u,v)}
\begin{pmatrix}
1+u[H_{2,1}(v)-H_{2,2}(v)] \\ 1-u[H_{1,1}(v)-H_{1,2}(v)]
\end{pmatrix}.
$$

\subsection{Case $\ell=r=1$}

We have in this case $\cE=\cR=\{a\}$. The sojourn time of interest reads now
$$
\Tna=\sum_{\i=1}^n \ind_{\{S(\i)=a\}}.
$$
We are dealing with the local time of the one-dimensional random walk at~$a$.
We adapt the settings by putting
\begin{align*}
p^{(\i)}
&
=\Pa\{S(\i)=a\}=\Po\{S(\i)=0\},
\\
q^{(\j)}
&
=\Pa\{\ta=\j\}=\Po\{\tau_{_0}=\j\},
\\
r^{(k,n)}
&
=\Pa\{\Tna=k\}=\Po\{\Tno=k\},
\end{align*}
and
$$
G(u)=\sum_{\i=0}^{\infty} p^{(\i)} u^{\i},\quad
H(u)=\sum_{\j=1}^{\infty} q^{(\j)} u^{\j}=\Eo\!\left(u^{\tau_{_0}}\right).
$$
We explicitly have
$$
p^{(\i)}=\binom{\i}{\i/2} (pq)^{\i/2},\quad
q^{(\j)}=\frac{1}{\j-1}\binom{\j}{\j/2} (pq)^{\j/2},
$$
$$
G(u)=\frac{1}{A(u)},\quad H(u)=1-\frac{1}{G(u)}=1-A(u).
$$
By~\refp{law-TnR}, we have
$$
\sum_{n=k}^{\infty} r^{(k,n)}v^n=\frac{A(v)}{1-v} \,[1-A(v)]^k.
$$
In this case, we can explicitly invert the previous generating function.
Indeed, by invoking~\refp{sum}, we get
\begin{align*}
\sum_{n=k}^{\infty} r^{(k,n)}v^n
&
=(1-4pqv^2)\left(\sum_{l=0}^{\infty}v^l\right)
\sum_{m=0}^{\infty} 2^k\binom{2m+k}{m}(pqv^2)^{m+k}
\\
&
=2^k\left(1+v+(1-4pq)\sum_{l=2}^{\infty}v^l\right)
\sum_{m=k}^{\infty} \binom{2m-k}{m}(pq)^m v^{2m}.
\end{align*}
The foregoing double sum can be easily transformed as follows:
$$
\sum_{l=2}^{\infty}v^l\sum_{m=k}^{\infty} \binom{2m-k}{m}(pq)^m v^{2m}
=\sum_{n=2k+2}^{\infty}
\Bigg(\sum_{m=k}^{[n/2]-1} \binom{2m-k}{m}(pq)^m\Bigg) v^n
$$
and then
\begin{align*}
\sum_{n=k}^{\infty} r^{(k,n)}v^n
&
=2^k\Bigg[\sum_{m=k}^{\infty} \binom{2m-k}{m}(pq)^m v^{2m}
+\sum_{m=k}^{\infty} \binom{2m-k}{m}(pq)^m v^{2m+1}
\\
&
\hphantom{=}+(1-4pq)\sum_{n=2k+2}^{\infty}
\Bigg(\sum_{m=k}^{[n/2]-1} \binom{2m-k}{m}(pq)^m\Bigg) v^n\Bigg]\!.
\end{align*}
Hence, by identifying the terms in~$v^n$, we find that for even~$n$
\begin{equation}\label{rkn-inter}
r^{(k,n)}=2^k \Bigg[\binom{n-k}{n/2}(pq)^{n/2} +(1-4pq)
\sum_{m=k}^{n/2-1} \binom{2m-k}{m}(pq)^m\Bigg]
\end{equation}
and for odd~$n$, observing that $r^{(k,n)}=r^{(k,n-1)}$ since the last
step of the random walk cannot vanish in this case,
$$
r^{(k,n)}=2^k \Bigg[\binom{n-k-1}{(n-1)/2}(pq)^{(n-1)/2} +(1-4pq)
\sum_{m=k}^{(n-3)/2} \binom{2m-k}{m}(pq)^m\Bigg]\!.
$$
In the two above formulas we adopted the convention that
$\sum_{m=k}^l=0$ when $k>l$.
At this stage, we see that in the particular case of the symmetric
random walk (corresponding to $p=q=\frac12$),
the foregoing formulas simply reduce to
$$
r^{(k,n)}=\begin{cases}
\displaystyle\frac{\binom{n-k}{n/2}}{2^{n-k}} & \mbox{if $n$ is even},
\\[2ex]
\displaystyle\frac{\binom{n-k-1}{(n-1)/2}}{2^{n-k-1}} & \mbox{if $n$ is odd}.
\end{cases}
$$
Going back to~\refp{rkn-inter}, we get for even~$n$:
\begin{align*}
r^{(k,n)}
&
=2^k \Bigg[\binom{n-k}{n/2}(pq)^{n/2}
+\sum_{m=k}^{n/2-1} \binom{2m-k}{m}(pq)^m
-4\sum_{m=k}^{n/2-1} \binom{2m-k}{m}(pq)^{m+1}\Bigg]
\\
&
=2^k \Bigg[\binom{n-k}{n/2}(pq)^{n/2}
+ \sum_{m=k}^{n/2-1} \binom{2m-k}{m}(pq)^m
-4\sum_{m=k+1}^{n/2} \binom{2m-k-2}{m-1}(pq)^m\Bigg]
\\
&
=2^k \Bigg[(pq)^k+\sum_{m=k+1}^{n/2} \left[\binom{2m-k}{m}-4\binom{2m-k-2}{m-1}\right](pq)^m\Bigg]
\\
&
=2^k \Bigg[(pq)^k+\sum_{m=k+1}^{n/2} \frac{(2m-k-2)!}{m!(m-k)!}\,(k^2+k-2m)(pq)^m\Bigg]
\\
&
=(2pq)^k \Bigg[1+\sum_{m=1}^{n/2-k} \frac{(2m+k-2)!}{m!(m+k)!}\,(k^2-k-2m)(pq)^m\Bigg].
\end{align*}
Consequently, we have obtained the following result.
%
\begin{theo}
The probability distribution of~$\Tno$ is given, for even~$n$, by
$$
\Po\{\Tno=k\}=(2pq)^k \Bigg[1+\sum_{m=1}^{n/2-k} \frac{(2m+k-2)!}{m!(m+k)!}\,(k^2-k-2m)(pq)^m\Bigg]
$$
and for odd~$n$ by $\Po\{\Tno=k\}=\Po\{\Tnuno=k\}$.
When $p=q=\frac12$ (case of the symmetric random walk), this distribution
reduces, for even~$n$, to
$$
\Po\{\Tno=k\}=\frac{\binom{n-k}{n/2}}{2^{n-k}}.
$$
\end{theo}
%
In particular, we have
$$
r^{0,n}=1-\sum_{m=1}^{n/2} \frac{\binom{2m}{m}}{2m-1}\,(pq)^m.
$$
This quantity, which represents the probability $\Po\{\Tno=0\}$ is nothing
but the probability $\Po\{\tzero>n\}$.

\section{Further investigations}

In a more realistic model, the plasmic membrane should be viewed as a closed curve
and its roaming constituents should be modeled for instance by random walks on
the finite torus $\Z/N\Z=\{0,1,\dots,N-1\}$ for a fixed possibly large integer~$N$
(with the usual rule $N\equiv 0$ $[\mathrm{mod}\,N]$) or by Brownian motions
on the continuous torus $\R/\Z$.

Moreover, the rafts could be built by aggregating several receptors, that is
by choosing several sequences of successive receptors:
$\{a_1,a_1+1,\dots,a_1+\l_1-1\}$,
$\{a_2,a_2+1,\dots,a_2+\l_2-1\},\dots,\{a_r,a_r+1,\dots,a_r+\l_r-1\}$
with $a_1+\l_1\le a_2$, $a_2+\l_2\le a_3,\dots,a_{r-1}+\l_{r-1}\le a_r$,
where the $\l_1,\dots,\l_r$ are the length of the rafts. The effect of this
kind of non-homogeneous repartition of receptors on the binding distribution
could provide a possible functional property of rafts. In addition, the
localizations and the lengths of the rafts could be random and then we should
also consider that the numbers $a_1,\dots,a_r$ and $\l_1,\dots,\l_r$ are random variables.

So far we concentrate on diffusion-limited reaction neglecting the binding
duration. But in the reality, ligand-receptor binding induces some delay in the mechanism,
that is when a ligand meets a receptor, they bind during a certain (possibly
random) amount of time before unbinding. So, we should also introduce convenient
delayed random walk or Brownian motion for recreating the real biological process.

In Subsection~\ref{subsect-walk}, we give some information about a possible
random walk model on the torus $\Z/N\Z$ and in Subsection~\ref{subsect-BM}, we address
possible continuous models involving Brownian motion on the line $\R$ or
on the torus $\R/\Z$.

\subsection{A random walk model on $\Z/N\Z$}\label{subsect-walk}

For building a  model on $\Z/N\Z$ with a deterministic set of receptors $\cR'=\{a_1,\dots,a_r\}$,
we introduce an $\ell$--dimensional random walk
$(\bS'(\i))_{\i\in\N}$ on the finite set $(\Z/N\Z)^{\ell}$. Set
$$
\cE'=\bigcup_{j=1}^{\ell} \left[(\Z/N\Z)^{j-1}\times\cR' \times(\Z/N\Z)^{\ell-j}\right].
$$
This set is finite: $\#\cE'=\ell r N^{\ell-1}$, and the sojourn time of interest writes here
$$
\Tnprime=\sum_{\i=1}^n \ind_{\{\bS'(\i)\in \cE'\}}.
$$
For tackling the computation of the probability distribution of $\Tnprime$,
we observe that we can pass from the one--dimensional random walk
$(S'(\i))_{\i\in\N}$ on $\Z/N\Z$ to a walk on~$\Z$ in the following manner:
set, for $\i\in\N$,
$$
S(\i)=S'(\i)+\a_{\i}N
$$
where $\a_{\i}$ is the number of upcrossings of level~$0$ (i.e. times~$\s$
such that $S'(\s-1)$ $=N-1,S'(\s)=0,S'(\s+1)=1$) minus that of downcrossing of $0$
(i.e. times~$\s$ such that $S'(\s-1)=1,S'(\s)=0,S'(\s+1)=N-1$) up to time~$\i$.
Then, $(S(\i))_{\i\in\N}$ is a random walk on $\Z$.
Conversely, $S'(\i)=S(\i)\!\mod N=S(\i) -N[S(\i)/N]$.
Notice that in our new setting, we mark with primes all quantities related to
the walks on $\Z/N\Z$ and $(\Z/N\Z)^{\ell}$, while those related
to the walks on $\Z$ and $\Z^{\ell}$ are not marked with any prime.

In this correspondence, the set of receptors $\cR'$ in~$\Z/N\Z$
becomes an infinite set $\cR=\cR'+N\Z=\{a_1+mN,\dots,a_r+mN;m\in\Z\}$
of receptors in~$\Z$ and the subset $\cE'$ of $(\Z/N\Z)^{\ell}$ becomes
the subset of~$\Z^{\ell}$
\begin{equation}\label{newset}
\cE=\bigcup_{j=1}^{\ell} \left(\Z^{j-1}\times\cR \times \Z^{\ell-j}\right)\!.
\end{equation}
Of course,
$$
\Tnprime=\sum_{\i=1}^n \ind_{\{\bS(\i)\in \cE\}}.
$$
The analysis done in the case of the walk on $\Z^{\ell}$ (associated with the
set~\refp{newset} now) can be carried out for the walk on $(\Z/N\Z)^{\ell}$ exactly
in the same way \textit{mutatis mutandis}:
the generating matrix~$\bK'(u,v)$ of the family of numbers
$\bPx\{\Tnprime=k\}$, $\bx\in\cE'$, $k,n\in\N$, is given by
$$
\bK'(u,v)=\frac{1}{1-v} \,[\bIEprime-u\bH'(v)]^{-1}[\bIEprime-\bH'(v)]\unEprime
$$
where $\bIEprime$ is the identity matrix on $\cE'$:
$\bIEprime=(\d_{\bx,\by})_{_{\scriptstyle \bx,\by\in\cE'}}$,
$\unEprime$ is the column-matrix consisting of~$1$,
that is $(1)_{_{\scriptstyle \bx\in\cE'}}$, and
$$
\bH'(v)=\bIEprime-\bG'(v)^{-1}.
$$
In the foregoing definition of $\bH'(v)$, $\bG'(v)$ is the matrix generating function
of the probabilities $\bPx'\{\bS'(\i)=\by\}$, $\bx,\by\in\cE'$ which are given by
$$
\bPx'\{\bS'(\i)=\by\}=\prod_{j=1}^{\ell} \Pxj'\{S'(\i)=y_j\}.
$$
The probabilities related to the case $\ell=1$ can be evaluated as follows.
For $x,y\in\Z/N\Z$,
\begin{align*}
\Px'\{S'(\i)=y\}
&
=\sum_{k\in\Z} \Px\{S(\i)=y+kN\}
\\
&
=\sum_{k\in\Z:\atop (x-y-\i)/N\le k\le (x-y+\i)/N} \binom{\i}{\frac{\i+y-x+kN}{2}} p^{(\i+y-x+kN)/2} q^{(\i+x-y-kN)/2}.
\end{align*}

Let us have a look to the case $\ell=r=1$. For this we introduce the following
families of hitting times related to $(S(\i))_{\i\in\N}$: for $a,b,c\in\Z$ such that $a<b<c$,
\begin{align*}
\tab & =\min\{\i\ge 1: S(\i)\in\{a,b\}\},
\\
\tabc & =\min\{\i\ge 1: S(\i)\in\{a,b,c\}\},
\end{align*}
and that related to $(S'(\i))_{\i\in\N}$: for $a\in\Z/N\Z$,
$$
\tau_a'=\min\{\i\ge 1: S'(\i)=a\}.
$$
It is clear that, if the starting point of $(S'(\i))_{\i\in\N}$ is $a$,
$$
\tau_a'=\tau_{\raisebox{-0.3ex}{$\scriptstyle a-N,a,a+N$}}.
$$
Because of this relationship, the probability distribution of time $\tau_a'$
can be explicitly written out. In the same way as~\refp{fg-tau-abc}, we have
$$
H'(v)=\E_a'\big(v^{\tau_a'}\big)
=\E_{\scriptscriptstyle 0}' \big(v^{\tau_{\raisebox{-0.1ex}{$\scriptscriptstyle 0$}}'}\big)
=\Eo\!\left(v^{\tau_{\raisebox{-0.1ex}{$\scriptscriptstyle -N,0,N$}}}\right)
=pv\,\E_{_1}\!\left(v^{\tau_{\raisebox{-0.1ex}{$\scriptscriptstyle 0,N$}}}\right)
+qv\,\E_{_{-1}}\!\left(v^{\tau_{\raisebox{-0.1ex}{$\scriptscriptstyle -N,0$}}}\right)\!.
$$
Invoking~\refp{fg-tau-ab-bis} and the fact that $B^+(v)B^-(v)=q/p$, we obtain
\begin{align*}
\E_{_1}\!\left(v^{\tau_{\raisebox{-0.1ex}{$\scriptscriptstyle 0,N$}}}\right)
&
=\frac{(1-B^-(v)^N)B^+(v)-(1-B^+(v)^N)B^-(v)}{B^+(v)^N-B^-(v)^N}
\\
&
=\frac{B^+(v)-B^-(v)}{B^+(v)^N-B^-(v)^N}
+\frac qp\frac{B^+(v)^{N-1}-B^-(v)^{N-1}}{B^+(v)^N-B^-(v)^N},
\end{align*}
\begin{align*}
\E_{_{-1}}\!\left(v^{\tau_{\raisebox{-0.1ex}{$\scriptscriptstyle 0,N$}}}\right)
&
=\frac{(1-B^-(v)^N)B^+(v)^{N-1}-(1-B^+(v)^N)B^-(v)^{N-1}}{B^+(v)^N-B^-(v)^N}
\\
&
=\frac{B^+(v)^{N-1}-B^-(v)^{N-1}}{B^+(v)^N-B^-(v)^N}
+\left(\frac qp\right)^{\!N-1} \frac{B^+(v)-B^-(v)}{B^+(v)^N-B^-(v)^N},
\end{align*}
which entails
$$
H'(v)=\bigg[\!\left(\frac qp\right)^{\!N}\!+1\bigg]pv\,\frac{B^+(v)-B^-(v)}{B^+(v)^N-B^-(v)^N}
+2qv\,\frac{B^+(v)^{N-1}-B^-(v)^{N-1}}{B^+(v)^N-B^-(v)^N}.
$$
Finally, the probability distribution of the sojourn time $\Tnoprime$ is
characterized by the generating function, as in~\refp{law-TnR},
$$
\sum_{n=k}^{\infty} \P_{\scriptscriptstyle 0}'\{\Tnoprime=k\} v^n
=\frac{1}{1-v} \,[H'(v)]^k [1-H'(v)].
$$
It seems difficult to extract the probabilities
$\P_{\scriptscriptstyle 0}'\{\Tnoprime=k\}$, $1\le k\le n$,
from this last identity.


\subsection{Brownian models on $\R$ and $\R/\Z$}\label{subsect-BM}

At last, we evoke the continuous counterpart to our model in the one-dimensional
case. Analogous continuous models hinge on linear Brownian motion (that is the
membrane is viewed as the real line $\R$) and circular Brownian motion (in this
case, the membrane is modeled as the torus $\R/\Z$).
In both cases, the rafts (clustered receptors) are viewed as intervals $[a_1,b_1],
\dots,[a_r,b_r]$ of $\R$ or $\R/\Z$, and ligands move like independent Brownian motions
$(B_1(s))_{s\ge 0},\dots,(B_{\ell}(s))_{s\ge 0}$ on $\R$ or $\R/\Z$. Set
$$
\mathbf{B}(s)=(B_1(s),\dots,B_{\ell}(s))\mbox{ for any $s\ge 0$},
$$
and
$$
\cR=\bigcup_{i=1}^r [a_i,b_i],\,\cE=\bigcup_{j=1}^{\ell} \left[\R^{j-1}\times\cR \times\R^{\ell-j}\right]
\mbox{ or }\bigcup_{j=1}^{\ell} \left[(\R/\Z)^{j-1}\times\cR \times(\R/\Z)^{\ell-j}\right]\!.
$$
The process $(\mathbf{B}(s))_{s\ge 0}$ is a Brownian motion on $\R^{\ell}$ or
$(\R/\Z)^{\ell}$. The time that ligands bind with rafts up to a fixed time~$t$ is
given by the sojourn time of $(\mathbf{B}(s))_{s\ge 0}$ in $\cE$:
$$
\Tt=\int_0^t \ind_{\{\mathbf{B}(s)\in\cE\}}\,\mathrm{d}s.
$$
In the case of linear Brownian motion ($\ell=1$), the computation of the probability
distribution of $\Tt$ will be the object of a forthcoming work (\cite{lachal}).

\appendix\section*{Appendix}
\section{Distribution of the first hitting time for the random walk}\label{appendix1}

The generating function~\refp{solutionH} writes
$$
H_{x,a,b}^+(u)=(2pu)^{b-x} \frac{[1+A(u)]^{x-a}-[1-A(u)]^{x-a}}{[1+A(u)]^{b-a}-[1-A(u)]^{b-a}}
$$
Let us introduce the rational fraction defined for $x,y\in\N$ such that $x<y$, by
$$
f_{x,y}(\z)=\frac{(1+\z)^x-(1-\z)^x}{(1+\z)^y-(1-\z)^y}
=\frac{\sum_{m=0}^{[(x-1)/2]} \binom{x}{2m+1} \z^{2m}}
{\sum_{m=0}^{[(y-1)/2]} \binom{y}{2m+1} \z^{2m}}.
$$
We want to decompose this fraction into partial fractions. The poles of $f_{x,y}$
are the roots of the polynomial $\sum_{m=0}^{[(y-1)/2]} \binom{y}{2m+1} \z^{2m}$
which are those of $(1+\z)^y-(1-\z)^y$ except for $0$. It is easily seen that
they consist of the family of numbers $\z_l=\frac{e^{2il\pi/y}-1}{e^{2il\pi/y}+1}=i\tan(l\pi/y)$
for $1\le l\le y-1$ such that $l\ne y/2$. Moreover, the function $\z\mapsto f_{x,y}(\z)$ is even,
so we decompose it as
$$
f_{x,y}(\z)=\sum_{l\in\N:\atop 1\le l<y/2} \frac{\a_l}{\z^2-\z_l^2}
$$
with $\a_l=\lim_{\z\to\z_l} (\z^2-\z_l^2) f(\z)$. We have
\begin{align*}
\a_l
&
=\frac{(1+\z_l)^x-(1-\z_l)^x}{\lim_{\z\to\z_l} \frac{(1+\z)^y-(1-\z)^y}{\z^2-\z_l^2}}
=\frac{2\z_l}{y}\,\frac{(1+\z_l)^x-(1-\z_l)^x}{(1+\z_l)^{y-1}+(1-\z_l)^{y-1}}
\\
&
=2(-1)^{l-1}\cos^{y-x-3}\!\left(\frac{l\pi}{y}\right)
\sin\!\left(\frac{l\pi}{y}\right) \sin\!\left(\frac{lx\pi}{y}\right)\!.
\end{align*}
Therefore,
$$
f_{x,y}(\z)=2\sum_{l\in\N:\atop 1\le l<y/2} (-1)^{l-1}\frac{\cos^{y-x-3}\!\big(\frac{l\pi}{y}\big)
\sin\!\big(\frac{l\pi}{y}\big)\sin\!\big(\frac{lx\pi}{y}\big)}{\z^2+\tan^2\!\big(\frac{l\pi}{y}\big)}.
$$
Let us apply this formula to $H_{x,a,b}^+(u)=(2pu)^{b-x} f_{x-a,b-a}(A(u))$.
Since
$$
A(u)^2+\tan^2\!\left(\frac{l\pi}{b-a}\right)
=\frac{1}{\cos^2\!\big(\frac{l\pi}{b-a}\big)}
\left[1-4pq\!\left(\cos^2\!\left(\frac{l\pi}{b-a}\right)\!\right)u^2\right]\!,
$$
we get that
$$
H_{x,a,b}^+(u)=2(2pu)^{b-x}\sum_{l\in\N:\atop 1\le l<(b-a)/2}
\frac{\cos^{b-x-1}\!\big(\frac{l\pi}{b-a}\big)
\sin\!\big(\frac{l\pi}{b-a}\big)\sin\!\big(\frac{l(b-x)\pi}{b-a}\big)}{
1-4pq\big[\cos^2\!\big(\frac{l\pi}{b-a}\big)\big]u^2}.
$$
Using the expansion
\begin{align*}
\frac{1}{1-4pq\,\big[\!\cos^2\!\big(\frac{l\pi}{b-a}\big)\big]u^2}
=\sum_{\j=0}^{\infty} \left[4pq\cos^2\!\left(\frac{l\pi}{b-a}\right)\!\right]^{\j} u^{2\j},
\end{align*}
we find
\begin{align*}
H_{x,a,b}^+(u)
&
=\sum_{\j=0}^{\infty}2(2p)^{b-x}
\Bigg[\sum_{l\in\N:\atop 1\le l<(b-a)/2} \cos^{b-x-1}\!\left(\frac{l\pi}{b-a}\right)
\sin\!\left(\frac{l\pi}{b-a}\right)
\\
&
\hphantom{=\;} \sin\!\left(\frac{l(b-x)\pi}{b-a}\right)\!
\left(4pq\cos^2\!\left(\frac{l\pi}{b-a}\right)\!\right)^{\j}\Bigg] u^{2\j+b-x}.
\end{align*}
Performing the change of index $\j\mapsto(\j+x-b)/2$ into the above sum,
we obtain
\begin{align*}
H_{x,a,b}^+(u)
&
=\sum_{\j=b-x}^{\infty} 2\left(\frac pq\right)^{\!(x-b)/2}
\Bigg[\sum_{l\in\N:\atop 1\le l<(b-a)/2} \cos^{b-x-1}\!\left(\frac{l\pi}{b-a}\right)
\sin\!\left(\frac{l\pi}{b-a}\right)
\\
&
\hphantom{=\;}\sin\!\left(\frac{l(b-x)\pi}{b-a}\right)\!
\left(2\sqrt{pq}\cos\!\left(\frac{l\pi}{b-a}\right)\!\right)^{\j}\Bigg] u^{\j}
\end{align*}
from which we finally derive the expression~\refp{distribution-tau-ab} of the
probability $q_{x,a,b}^{(\j)+}$. In a quite analogous way, we can deduce the
probability $q_{x,a,b}^{(\j)-}$ from $H_{x,a,b}^-(u)=(2qu)^{x-a} f_{b-x,b-a}(A(u))$.

\section{Distribution of the stopped random walk}\label{appendix2}

As mentioned in Section~\ref{background2}, the stopping-probabilities can be
evaluated by invoking the famous reflection principle. We provide the details
in this Appendix.

\vspace{.5\baselineskip}
\noindent\textbf{\textsl{One-sided threshold}}
\vspace{.5\baselineskip}

Our aim is to compute the probability of the random walk stopped
when hitting the threshold $a$: $\Px\{S(\j)=y,\j\le \ta\}$. We begin by
first evaluating $\Px\{S(\j)=y,\ta<\j\}$.
We have, for $x,y\in\Z$ such that $x,y<a$ or $x,y>a$,
$$
\Px\{S(\j)=y,\ta<\j\}=\sum_{\i=1}^{\j-1}\Px\{\ta=\i\} \Pa\{S(\j-\i)=y\}.
$$
The reflection principle stipulates that
$$
\Pa\{S(\j-\i)=y\}=\left(\frac pq \right)^{\!y-a}\Pa\{S(\j-\i)=2a-y\}.
$$
Hence,
\begin{align*}
\Px\{S(\j)=y,\ta<\j\}
&
=\sum_{\i=1}^{\j-1}\Px\{\ta=\i\} \left(\frac pq \right)^{\!y-a}\Pa\{S(\j-\i)=2a-y\}
\\
&
=\left(\frac pq \right)^{\!y-a}\Px\{S(\j)=2a-y,\ta<\j\}
\\
&
=\left(\frac pq \right)^{\!y-a}\Px\{S(\j)=2a-y\}.
\end{align*}
In the last equality, we used the fact that $x$ and $2a-y$ are around $a$
and then the condition $\ta<\j$ is redundant. Therefore,
$$
\Px\{S(\j)=y,\j\le \ta\}
=\Px\{S(\j)=y\}-\left(\frac pq \right)^{\!y-a}\Px\{S(\j)=2a-y\}
$$
which is nothing but~\refp{killed-a}.

\vspace{.5\baselineskip}
\noindent\textbf{\textsl{Two-sided threshold}}
\vspace{.5\baselineskip}

Here, we aim to compute the probability of the random walk stopped
when hitting one of the two thresholds $a$ and $b$: $\Px\{S(\j)=y,\j\le \tab\}$
under the assumption $x\in(a,b)$.
For this, we introduce the successive hitting times of levels $a$ and $b$
(with the usual convention $\min\emptyset=+\infty$):
$$
\begin{array}{lll}
\so=\ta,& \s_{2l+1}=\min\{\j>\s_{2l}:S(\j)=b\},
&\s_{2l+2}=\min\{\j>\s_{2l+1}:S(\j)=a\},
\\
\vso=\tb,& \vs_{2l+1}=\min\{\j>\vs_{2l}:S(\j)=a\},
&\vs_{2l+2}=\min\{\j>\vs_{2l+1}:S(\j)=b\}.
\end{array}
$$
Notice that if $\ta<\tb$, then $\s_{l+1}=\vs_l$
and if $\ta>\tb$, then $\vs_{l+1}=\s_l$. In both cases, we see that
$$
\min(\s_{l},\vs_{l})=\begin{cases}\s_l &\mbox{if } \ta<\tb\\
\vs_l &\mbox{if } \ta>\tb\end{cases},\quad
\max(\s_{l},\vs_{l})=\begin{cases}\s_{l+1} &\mbox{if } \ta<\tb\\
\vs_{l+1} &\mbox{if } \ta>\tb\end{cases}
$$
and that
$$
\max(\s_{l},\vs_{l})=\min(\s_{l+1},\vs_{l+1}).
$$
With these properties at hands, we write the event $\{\tab<\j\}$ as follows:
\begin{align*}
\{\tab<\j\}
&
=\{\ta<\tb,\ta<\j\}\cup \{\ta>\tb,\tb<\j\}
\\
&
=\left(\bigcup_{l\in\N} \{\ta<\tb,\s_{2l}<\j\le\s_{2l+2}\}\right)
\cup\left(\bigcup_{l\in\N} \{\ta>\tb,\vs_{2l}<\j\le\vs_{2l+2}\}\right)
\\
&
=\left(\bigcup_{l\in\N} \{\ta<\tb,\s_{2l}<\j\le\vs_{2l+1}\}\right)
\cup\left(\bigcup_{l\in\N} \{\ta>\tb,\vs_{2l}<\j\le\s_{2l+1}\}\right)
\\
&
=\bigcup_{l\in\N} \{\min(\s_{2l},\vs_{2l})<\j\le\max(\s_{2l+1},\vs_{2l+1})\}.
\end{align*}
Now, the probability of the stopped random walk can be written as
\begin{align}
\lqn{\Px\{S(\j)=y,\j\le\tab\}}
&
=\Px\{S(\j)=y\}-\sum_{l=0}^{\infty} \Px\{S(\j)=y,
\min(\s_{2l},\vs_{2l})<\j\le\max(\s_{2l+1},\vs_{2l+1})\}.
\label{killed-ab-bis}
\end{align}
Let us evaluate the term lying within the sum~\refp{killed-ab-bis}:
\begin{align}
\lqn{\Px\{S(\j)=y,\min(\s_{2l},\vs_{2l})<\j\le\max(\s_{2l+1},\vs_{2l+1})\}}
&
=\Px\{S(\j)=y,\min(\s_{2l},\vs_{2l})<\j\}
-\Px\{S(\j)=y,\max(\s_{2l+1},\vs_{2l+1})<\j\}
\nonumber\\
&
=\Px\{(S(\j)=y,\s_{2l}<\j)\cup(S(\j)=y,\vs_{2l}<\j)\}
\nonumber\\
&
\hphantom{=\;}-\Px\{(S(\j)=y,\s_{2l+1}<\j)\cap(S(\j)=y,\vs_{2l+1}<\j)\}
\nonumber\\
&
=\Px\{S(\j)=y,\s_{2l}<\j\}+\Px\{S(\j)=y,\vs_{2l}<\j\}
\nonumber\\
&
\hphantom{=\;}-\Px\{(S(\j)=y,\s_{2l}<\j)\cap(S(\j)=y,\vs_{2l}<\j)\}
\nonumber\\
&
\hphantom{=\;}-\Px\{S(\j)=y,\s_{2l+1}<\j\}-\Px\{S(\j)=y,\vs_{2l+1}<\j\}
\nonumber\\
&
\hphantom{=\;}+\Px\{(S(\j)=y,\s_{2l+1}<\j)\cup(S(\j)=y,\vs_{2l+1}<\j)\}.
\label{sum1}
\end{align}
We observe that
\begin{align*}
\lqn{\Px\{(S(\j)=y,\s_{2l+1}<\j)\cup(S(\j)=y,\vs_{2l+1}<\j)\}}
&
=\Px\{S(\j)=y,\min(\s_{2l+1},\vs_{2l+1})<\j\}
\\
&
=\Px\{S(\j)=y,\max(\s_{2l},\vs_{2l})<\j\}
\\
&
=\Px\{(S(\j)=y,\s_{2l}<\j)\cap(S(\j)=y,\vs_{2l}<\j)\}.
\end{align*}
Thus, two terms of the sum~\refp{sum1} cancel and it remains
\begin{align}
\lqn{\Px\{S(\j)=y,\min(\s_{2l},\vs_{2l})<\j\le\max(\s_{2l+1},\vs_{2l+1})\}}
&
=\Px\{S(\j)=y,\s_{2l}<\j\}+\Px\{S(\j)=y,\vs_{2l}<\j\}
\nonumber\\
&
\hphantom{=\;}-\Px\{S(\j)=y,\s_{2l+1}<\j\}-\Px\{S(\j)=y,\vs_{2l+1}<\j\}.
\label{sum2}
\end{align}
By using the principle of reflection with respect to level $a$, we get
for the first term of the sum~\refp{sum2}:
\begin{align}
\Px\{S(\j)=y,\s_{2l}<\j\}
&
=\sum_{\i=1}^{\j-1} \Px\{\s_{2l}=\i\}\Pa\{S(\j-\i)=y\}
\nonumber\\
&
=\sum_{\i=1}^{\j-1}\Px\{\s_{2l}=\i\} \left(\frac pq \right)^{\!y-a}\Pa\{S(\j-\i)=2a-y\}
\nonumber\\
&
=\left(\frac pq \right)^{\!y-a}\Px\{S(\j)=2a-y,\s_{2l}<\j\}
\nonumber\\
&
=\left(\frac pq \right)^{\!y-a}\Px\{S(\j)=2a-y,\s_{2l-1}<\j\}.
\label{reflec1}
\end{align}
In the last step, we used the fact than, when $x,y\in(a,b)$, $2a-y<a$
and then the condition $S(\j)=2a-y,\s_{2l-1}<\j$ is enough to assure
$\s(2l)<\j$. Similarly, reflection with respect to level $b$ yields
\begin{align}
\lqn{\Px\{S(\j)=2a-y,\s_{2l-1}<\j\}}
=\left(\frac pq \right)^{\!2a-b-y}\Px\{S(\j)=y+2(b-a),\s_{2l-2}<\j\}.
\label{reflec2}
\end{align}
So, the following recursive relationship entails from~\refp{reflec1} and~\refp{reflec2}:
\begin{equation}\label{reflec3}
\Px\{S(\j)=y,\s_{2l}<\j\}=\left(\frac pq \right)^{\!a-b}\Px\{S(\j)=y+2(b-a),\s_{2l-2}<\j\}.
\end{equation}
By iterating~\refp{reflec3} with respect to index $l$, we plainly obtain
\begin{equation}\label{reflec4}
\Px\{S(\j)=y,\s_{2l}<\j\}=\left(\frac pq \right)^{\!l(a-b)}\Px\{S(\j)=y+2l(b-a),\ta<\j\}.
\end{equation}
A last application of the principle of reflection to~\refp{reflec4} supplies
\begin{equation}\label{reflec5}
\Px\{S(\j)=y,\s_{2l}<\j\}=\left(\frac pq \right)^{\!y-a+l(b-a)}\Px\{S(\j)=2a-2l(b-a)-y\}.
\end{equation}
Analogously, concerning the second term of the sum~\refp{sum2}, we
successively have
\begin{align}
\lqn{\Px\{S(\j)=y,\s_{2l+1}<\j\}}
&
=\left(\frac pq \right)^{\!l(b-a)}\Px\{S(\j)=y-2l(b-a),\s_1<\j\}
\nonumber\\
&
=\left(\frac pq \right)^{\!y-b-l(b-a)}\Px\{S(\j)=2b+2l(b-a)-y,\ta<\j\}
\nonumber\\
&
=\left(\frac pq \right)^{\!(l+1)(b-a)}\Px\{S(\j)=y-2(l+1)(b-a)\}.
\label{reflec6}
\end{align}
We obtain in a quite similar manner the two last terms of~\refp{sum2}:
\begin{align}
\lqn{\Px\{S(\j)=y,\vs_{2l}<\j\}}
&
=\left(\frac pq \right)^{\!y-b-l(b-a)}\Px\{S(\j)=2b+2l(b-a)-y\}
\nonumber\\
&
=\left(\frac pq \right)^{\!y-a-(l+1)(b-a)}\Px\{S(\j)=2a+2(l+1)(b-a)-y\},
\nonumber\\[-1ex]
\label{reflec7}\\[-1ex]
\lqn{\Px\{S(\j)=y,\vs_{2l+1}<\j\}}
&
=\left(\frac pq \right)^{\!-(l+1)(b-a)}\Px\{S(\j)=y+2(l+1)(b-a)\}.
\nonumber
\end{align}
Finally, by summing~\refp{reflec5}, \refp{reflec6} and \refp{reflec7}, we
deduce from~\refp{killed-ab-bis} and~\refp{sum2}:
\begin{align*}
\lqn{\Px\{S(\j)=y,\tab<\j\}}
&
=\sum_{l=0}^{\infty} \left[\left(\frac pq \right)^{\!y-a+l(b-a)} p_{x,2a-2l(b-a)-y}^{(\j)}
-\left(\frac pq \right)^{\!(l+1)(b-a)} p_{x,y-2(l+1)(b-a)}^{(\j)}\right.
\\
&
\hphantom{=\;}
\left.+\left(\frac pq \right)^{\!y-a-(l+1)(b-a)} p_{x,2a+2(l+1)(b-a)-y}^{(\j)}
-\left(\frac pq \right)^{\!-(l+1)(b-a)} p_{x,y+2(l+1)(b-a)}^{(\j)}\right]
\\
&
=\sum_{l=-\infty}^{\infty} \left(\frac pq \right)^{\!l(b-a)}
\left[\left(\frac pq \right)^{\!y-a} p_{x,2a-2l(b-a)-y}^{(\j)}
- p_{x,y-2l(b-a)}^{(\j)}\right]+\Px\{S(\j)=y\}
\end{align*}
from which~\refp{killed-ab} ensues.


\vspace{\baselineskip}
\noindent
\textsc{Acknowledgments.}
I thank H\'edi Soula for having addressed this problem to me and for his
help in writing the biological context in the introduction.


\begin{thebibliography}{99}

\bibitem{abram} Abramowitz, M.; Stegun, I. A.
Handbook of mathematical functions with formulas, graphs, and mathematical tables.
Dover Publications, 1972.

\bibitem{hedi} Car\'e, B. R.; Soula, H. A. Impact of receptors clustering on ligand binding.
BMC System Biology 2011, 5:48.

\bibitem{feller}
Feller, W. An introduction to probability theory and its applications.
Vol.~I. Third edition. John Wiley \& Sons, 1968.

\bibitem{lachal}
Lachal, A. Sojourn time in an union of intervals for diffusions. Submitted.

\bibitem{renyi} R\'enyi, A. Calcul des probabilit\'es. Dunod, 1966.

\bibitem{spitzer} Spitzer, F. Principles of random walk. Second edition.
Graduate Texts in Mathematics, Vol.~34.  Springer-Verlag, 1976.

\end{thebibliography}
\end{document}